\documentclass{amsart}
\usepackage{amssymb,amscd,amsthm,amsmath,amsfonts,eufrak}

\newcommand{\emp}{\emptyset}
\newcommand{\norm}[1]{\mbox{$\left\|#1\right\|$}}
\newcommand{\x}{\times}
\newcommand{\cs}{\mbox{$C^{*}$-algebra}}
\newcommand{\css}{\mbox{$C^{*}$-algebras}}

\newcommand{\C}{\mathbb{C}}

\newcommand{\R}{\mathbb{R}}

\newcommand{\ov}[1]{\mbox{$\overline{#1}$}}

\newcommand{\al}{\mbox{$\alpha$}}
\newcommand{\eps}{\mbox{$\epsilon$}}
\newcommand{\bt}{\mbox{$\beta$}}
\newcommand{\ga}{\mbox{$\gamma$}}
\newcommand{\Ga}{\mbox{$\Gamma$}}
\newcommand{\de}{\mbox{$\delta$}}
\newcommand{\De}{\mbox{$\Delta$}}

\newcommand{\la}{\mbox{$\lambda$}}

\newcommand{\Si}{\mbox{$\Sigma$}}
\newcommand{\si}{\mbox{$\sigma$}}

\newcommand{\mfA}{\mathfrak{A}}
\newcommand{\mfB}{\mathfrak{B}}

\newcommand{\mfH}{\mathfrak{H}}

\newcommand{\bgc}{\begin{center}}

\newcommand{\edc}{\end{center}}
\newcommand{\be}{\begin{enumerate}}
\newcommand{\ee}{\end{enumerate}}
\newcommand{\beq}{\begin{equation}}
\newcommand{\eeq}{\end{equation}}
\newcommand{\beqn}{\begin{eqnarray}}
\newcommand{\eeqn}{\end{eqnarray}}
\newcommand{\beqns}{\begin{eqnarray*}}
\newcommand{\eeqns}{\end{eqnarray*}}
\newcommand{\bq}{\begin{quote}}
\newcommand{\eq}{\end{quote}}

\newcommand{\bi}{\begin{itemize}}
\newcommand{\ei}{\end{itemize}}
\newcommand{\bd}{\begin{description}}
\newcommand{\ed}{\end{description}}

\newcommand{\lan}{\mbox{$\langle$}}
\newcommand{\ran}{\mbox{$\rangle$}}

\theoremstyle{plain}
\newtheorem{theorem}{Theorem}

\newtheorem{proposition}{Proposition}
\newtheorem{corollary}{Corollary}

\numberwithin{equation}{section}

\begin{document}
\title{The E-theoretic descent functor for groupoids}
\author{Alan L. T. Paterson}
\address{1001 Harlan Drive, Oxford, Mississippi 38655, USA}
\email{m1p1t@yahoo.com}
\date{April, 2007}
\begin{abstract}
The paper establishes, for a wide class of locally compact groupoids $\Ga$, the E-theoretic descent functor at the 
$\cs$ level, in a way parallel to that established for locally compact groups by Guentner, Higson and  Trout.  The second section shows that $\Ga$-actions on a $C_{0}(X)$-algebra $B$, where $X$ is the unit space 
of $\Ga$, can be usefully formulated in terms of an action on the associated bundle $B^{\sharp}$.  The third section shows that the functor $B\to C^{*}(\Ga,B)$ is continuous and exact, and uses the disintegration theory of J. Renault.  The last section establishes the existence of the descent functor under a very mild condition on $\Ga$, the main technical difficulty involved being that of finding a $\Ga$-algebra that plays the role of $C_{b}(T,B)^{cont}$ in the group case.  
\end{abstract}

\maketitle

\section{Introduction}

In a number of situations, in particular for the assembly map, the Baum-Connes conjecture and index theory 
(\cite[Theorem 3.4]{Kasp2},\cite{Kasp1,Connesbook,Tunov,Tumoy,Tupre} and many others) the descent homomorphism 
$j_{G}:KK_{G}(A,B)\to KK(C^{*}(G,A),C^{*}(G,B))$, where $G$ is a locally compact group and $A, B$ are $G$-$\css$, is of great importance.   (There is a corresponding result for the reduced crossed product algebras.)  In noncommutative geometry, classical group symmetry does not suffice, and one requires smooth groupoids in place of Lie groups 
(\cite{Cointeg,Connesbook}), so that it is important to have available constructions, such as that which gives the descent homomorphism, for groupoid, rather than just for group, actions.  To this end, the work of LeGall (\cite[7.2]{LeGall}) shows the existence of the $KK_{\Ga}$-descent homomorphism for $\Ga$-algebras, where $\Ga$ is a locally compact, $\si$-compact 
Hausdorff groupoid with left Haar system. 

A similar issue arises when we consider E-theory rather than KK-theory.  Non-equivariant E-theory is developed in 
\cite{CoHig,Connesbook,Blackadar}.  Guentner, Higson and Trout gave a definitive account of 
group equivariant E-theory in their memoir \cite{GHT}.  In particular (\cite[pp.47, 60ff.]{GHT}), they established the group equivariant E-theoretic descent functor and used it in their definition of the E-theoretic assembly map.  Another situation where the E-theoretic descent homomorphism is required is in the Bott periodicity
theorem for infinite dimensional Euclidean space which was established by Higson, Kasparov and Trout (\cite{HKT}), with its applications to the equivariant topological index and the Novikov higher signature conjecture.  The descent homomorphism associates to an equivariant asymptotic morphism from $A$ to $B$ a canonical homomorphism from $K_{G}(A)$ to $K_{G}(B)$, and this is how it is used in \cite{HKT}.   The present paper studies the descent homomorphism in the much more general situation involving groupoids rather than groups.   By modifying the method of \cite{GHT}, we prove the existence of the groupoid descent homomorphism at the $\cs$ level for a very wide class of groupoids.  

We start by reformulating the concept of a $\Ga$-action on a $C_{0}(X)$-algebra.  A $\Ga$-$\cs$ is defined in the literature (\cite{LeGall,Popescu}) as follows.   Let $X$ be the unit space of $\Ga$ and $A$ be a $C_{0}(X)$-algebra.  Roughly, the latter means that $A$ can be regarded as the 
$C_{0}(X)$-algebra of continuous sections vanishing at infinity of a $C^{*}$-bundle $A^{\sharp}$ of $\css$ $A_{x}$ 
$(x\in X)$.  One pulls back $A$ to $\Ga$ using the range and source maps $r, s$ to obtain $C_{0}(\Ga)$-algebras 
$r^{*}A, s^{*}A$.  An action of $\Ga$ on $A$ is just a 
$C_{0}(\Ga)$-isomorphism $\al:s^{*}A\to r^{*}A$ for which the maps $\al_{\ga}:A_{s(\ga)}\to A_{r(\ga)}$ compose in accordance with the rules for groupoid multiplication.  It is desirable to have available an equivalent definition for a
$\Ga$-action along the lines of an action (in the usual sense) of a group on a $\cs$: in groupoid terms, this should involve a continuity condition for the map 
$\ga\to \al_{\ga}$ on the $C^{*}$-bundle $A^{\sharp}$.  The specification of this continuity is very natural: we require that for each $a\in A$, the map $\ga\to \al_{\ga}(a_{s(\ga)})$ be continuous.  This definition is useful in a number of contexts, for example, in specifying the $\Ga$-algebra of continuous elements in a $\cs$ with an algebraic $\Ga$-action, and in working with the covariant algebra $C^{*}(\Ga,A)$.

The second section proves that the two definitions of $\Ga$-action on $A$ are equivalent.  We survey the theory of 
$C^{*}$-bundles, in particular, the topologizing of $A^{\sharp}$.   We require the well-known result, related to the Dauns-Hoffman theorem, that the ``Gelfand transform'' of $A$ is an isomorphism  onto $C_{0}(X,A^{\sharp})$.  A simple modification of a corresponding result by Dupr\'{e} and Gillette (\cite{DupreG}) gives this result and we sketch it for completeness.  (Another approach to this is given by Nilsen (\cite{Nilsen}).)  

Following the method of Guentner, Higson and Trout, we have to show that the functor $A\to C^{*}(\Ga,A)$ is continuous and exact.  This is proved in the third section.  The continuity of this functor is proved in a way similar to that of the group case, while for exactness, we give a groupoid version of the corresponding theorem of N. C. Phillips (\cite{PhillipsLNS}) for locally compact groups.  The proofs of these use the disintegration theorem of J. Renault.  

The fourth section establishes our version of the descent homomorphism.  The theory of groupoid equivariant asymptotic morphisms for $\Ga$-algebras is developed using the work of R. Popescu (\cite{Popescu}), to whom I am grateful for helpful comments.  In particular, in place of the non-equivariant $C_{b}(T,B)$, the $C_{0}(X)$-algebra 
$C_{b}^{X}(T,B)=C_{0}(X)C_{b}(T,B)$ is used.  A technical difficulty arises since (unlike the locally compact group case) there does not exist a natural algebraic $\Ga$-action on 
$C_{b}^{X}(T,B)$.  However, there is another natural bundle 
$C_{b}(T,B^{\sharp})=\sqcup C_{b}(T,B_{x})$ on which there is a simple algebraic
$\Ga$-action, derived from the given action on $B$, and a bundle map $R$ from $C_{b}^{X}(T,B)^{\sharp}$ to 
$C_{b}(T,B^{\sharp})$.  We show that if $\Ga$ has local $G$-sets - a very mild, ``transversal'' condition satisfied by most groupoids that arise in practice - then we {\em can} find a canonical $\Ga$-algebra $\mfB\subset C_{b}^{X}(T,B)$ determined by the action of $\Ga$ on $C_{b}(T,B^{\sharp})$.  The map $B\to \mfB$ is functorial and is the natural choice 
for defining the equivariant asymptotic algebra and the groupoid descent homomorphism.

\section{Groupoid $\css$}
 
Let $\Ga$ be a locally compact, second countable, Hausdorff groupoid with a left Haar system $\la$.  The unit space of 
$\Ga$ is denoted by $X$.  The range and source maps $r,s:\Ga\to X$ are given by: $r(\ga)=\ga\ga^{-1}$, 
$s(\ga)=\ga^{-1}\ga$.  We now review $\Ga$-spaces.  

A $\Ga$-space is a topological space $M$ with an open, onto, projection map $p:M\to X$, $M_{x}=p^{-1}(x)$, and a continuous product map (a {\em $\Ga$-action}) from \\
$\Ga\x_{s} M=\{(\ga,m): s(\ga)=p(m)\}\to M$, so that the usual groupoid algebra holds: in particular, if 
$s(\ga_{2})=p(m)$, then
$\ga_{1}(\ga_{2}m)=(\ga_{1}\ga_{2})m$ whenever $s(\ga_{1})=r(\ga_{2})$, and $p(\ga_{2} m)=r(\ga_{2})$. 
In the case where each $M_{x}$ is a $\cs$ and each of the maps 
$z\to \ga z$ is a $^{*}$-isomorphism from 
$M_{s(\ga)}$ onto $M_{r(\ga)}$, then we say that $M$ is a {\em $\Ga$-space of $\css$}.  We often write 
$\al_{\ga}(z)$ in place of $\ga z$.  For such an $M$, there is a natural groupoid $Iso(M)$ whose elements are the 
$^{*}$-isomorphisms $t$ from some $M_{x_{1}}$ onto $M_{x_{2}}$ and with unit space $X$.  Of course, 
$s(t)=x_{1}$ and $r(t)=x_{2}$ and the product is given by composition of maps.  Then saying that 
$\ga\to \al_{\ga}$ is a $\Ga$-action is equivalent to saying that the map is a groupoid homomorphism from $\Ga$ into
$Iso(M)$.   We call such an $M$ a {\em continuous} $\Ga$-space of $\css$ if the map 
$(\ga,z)\to \al_{\ga}(z)$ from  $\Ga\x_{s} M$ into $M$ is continuous.   

Let $A$ be a separable $\cs$.  We recall what it means for $A$ to be {\em $C_{0}(X)$-algebra}
(\cite{Kasp2,Blanchard92,Blanchard95,DupreG,LeGall,Nilsen}).   It means that
there is a homomorphism $\theta$ from $C_{0}(X)$ into the center $ZM(A)$ of the multiplier algebra  
$M(A)$ of $A$ such that $\theta(C_{0}(X))A=A$.  A $C_{0}(X)$-algebra $A$ determines a family of $\css$ $A_{x}$ ($x\in X$) where $A_{x}=A/(I_{x}A)$ with $I_{x}=\{f\in C_{0}(X): f(x)=0\}$.  If $J$ is a closed ideal of such an $A$, then the restriction map $f\to \theta(f)_{\mid J}$ makes $J$ also into a $C_{0}(X)$-algebra.  Also, $A/J$ is a $C_{0}(X)$-algebra in the obvious way, and $(A/J)_{x}=A_{x}/J_{x}$.

A {\em $C_{0}(X)$-morphism} from $A$ to $B$, where $A,B$ are $C_{0}(X)$-algebras, is defined to be a $^{*}$-homomorphism 
$T:A\to B$ which is also a $C_{0}(X)$-module map.  In that case, $T$ determines a 
$^{*}$-homomorphism $T_{x}:A_{x}\to B_{x}$ for each $x\in X$.  The following discussion is very close to, but not quite contained, in the book on Banach bundles by Dupr\'{e} and Gillette (\cite{DupreG}), and we will give a brief description of the modifications required.  Some of the details will be needed later and do not seem to appear in the literature.

Let $a_{x}$ be the image of $a\in A$ in $A_{x}$.  Let 
$A^{\sharp}=\sqcup A_{x}$ and $p:A^{\sharp}\to X$ be the map: $p(a_{x})=x$.  Then (cf. \cite[p.8]{DupreG}), $A^{\sharp}$ is a {\em $\cs$ family}: the map $p:A^{\sharp}\to X$ is surjective, and each fiber
$A_{x}=p^{-1}(x)$ is a $\cs$.  Let $\hat{a}$ be the section of $A^{\sharp}$ given by:
$\hat{a}(x)=a_{x}$.   By a $C^{*}$-family $E$ being a {\em $C^{*}$-bundle}, we mean 
(\cite[pp.6-9]{DupreG}) that $E$ is a topological space with $p$ open and continuous, that scalar multiplication, addition, multiplication and involution are continuous respectively from $\C\x E\to E$, from 
$E\x_{X} E\to E$, from $E\x_{X} E\to E$ and from 
$E\to E$, and the norm map $\norm{.}:E\to \R$ is upper semicontinuous and the following condition on the open sets for $E$
holds: if $W$ is open in 
$E$, $x\in X$ and the zero $0_{x}$ of $E_{x}$ belongs to $W$, then there exists an $\eps>0$ and an open neighborhood $U$ of $x$ such that 
\[   \{b\in p^{-1}(U): \norm{b}<\eps\}\subset W.                    \]
(Recall that a map $f:E\to \R$ is upper semicontinuous if, for each $e_{0}\in E$ and each $\eps>0$, there is an open neighborhood $U$ of $e_{0}$ in $X$ such that $f(e)<f(e_{0})+\eps$ for all $e\in U$.)

We now discuss briefly how, in a natural way, $E=A^{\sharp}$ is a $C^{*}$-bundle, and the map 
$a\to \hat{a}$ is a $^{*}$-isomorphism from $A$ onto $C_{0}(X,A^{\sharp})$ (a ``Gelfand'' theorem) 
(cf. \cite[2.1.3]{LeGall}, \cite{Nilsen}).   One proves first that each $\hat{a}$ is upper semicontinuous.  To this end, we modify the proof of the corresponding results 
(\cite[Proposition 2.1, Corollary 2.2]{DupreG}) which are proved in \cite{DupreG} for the completely regular, rather than locally compact Hausdorff, case.  The first of these results for our case can be stated as follows.  {\em For each $x\in X$, let 
$\mathcal{N}_{x}$ be the family of relatively compact, open subsets $V$ of $X$ containing $x$.  For each 
$V\in \mathcal{N}_{x}$, let $f_{V}:X\to [0,1]$ be continuous and such that it is $1$ on a neighborhood of $x$ in $V$ and vanishes outside $V$.  Then} 
\begin{equation}  \label{eq:afV}
\norm{\hat{a}(x)}=\inf\{\norm{f_{V}a}:V\in \mathcal{N}_{x}\}.    
\end{equation}
The upper semicontinuity of $\norm{\hat{a}}$ follow since if $\eps>0$ and $V$ is chosen so that 
$\norm{\hat{a}(x)}>\norm{f_{V}a}-\eps$, then for $y$ in a neighborhood of $x$, 
$\norm{\hat{a}(y)}\leq \norm{f_{V}a}<\norm{\hat{a}(x)}+\eps$.  (The equality (\ref{eq:afV}) is due to J. Varela.)

The first part of the proof of  (\ref{eq:afV}) shows that 
$\norm{\hat{a}(x)}\geq\inf\{\norm{f_{V}a}:V\in \mathcal{N}_{x}\}$.  This is the same as in the original Proposition 2.1.  For the reverse inequality, let $\eps>0$.  Since $A$ is a $C_{0}(X)$-algebra, $a=fb$ for some $f\in C_{0}(X), b\in A$, and 
using a bounded approximate identity in $C_{0}(X)$, there exists $F\in C_{0}(X)$ such that $0\leq F\leq 1$, $F(x)=1$ and 
$\norm{(1-F)a}<\eps$.  As $(F-f_{V})a\in I_{x}A$, and $a=f_{V}a + (1-F)a + (F-f_{V})a$,
\[  \norm{\hat{a}(x)}\leq \norm{f_{V}a + (1-F)a}\leq \norm{f_{V}a}+\eps  \]
and we obtain $\norm{\hat{a}(x)}\leq\inf\{\norm{f_{V}a}:V\in \mathcal{N}_{x}\}$.

Since we have (\ref{eq:afV}), the conditions of (\cite[Proposition 1.3]{DupreG}) 
(or of (\cite[Proposition 1.6]{Fell}) are satisfied, and $A^{\sharp}$ is a 
$C^{*}$-bundle.  A base for the topology on $A^{\sharp}$ (\cite[pp. 9-10, 16]{DupreG}) is given by sets of the form
\begin{equation}  \label{eq:aepsU} 
U(a,\eps)=\{b_{x}\in A_{x}: x\in U, \norm{b_{x}-a_{x}}<\eps\} 
\end{equation}
where $a\in A, \eps>0$ and $U$ is an open subset of $X$.  Further, a local base at $z\in A_{x_{0}}$ is given by neighborhoods of the form $U(a,\eps)$ where $a$ is any fixed element of $A$ for which $a_{x_{0}}=z$, and
$\hat{A}$ is a closed $^{*}$-subalgebra of $C_{0}(X,A^{\sharp})$.  To see that $\hat{A}=C_{0}(X,A^{\sharp})$, we just have to show (cf. (\cite[Proposition 2.3]{DupreG})) that $\hat{A}$ is dense in 
$C_{0}(X,A^{\sharp})$.  This follows by a simple partition of unity argument 
(\cite[Proposition 1.7]{Fell}, \cite[Lemma 5.3]{RaeWill}).
We then have the following theorem (\cite[Theor\`{e}me 2.1.1]{LeGall}).  It is also proved by Nilsen (\cite[Theorem 2.3]{Nilsen}) who derives the Dauns-Hoffman theorem (\cite[Theorem A.34]{RaeWill}) from it.  

\begin{theorem}   \label{th:asharp}
With the above topology on $A^{\sharp}$, $A^{\sharp}$ is a $C^{*}$-bundle over $X$.  Further, the relative topology on each $A_{x}$ is the norm topology.  Last, the map $a\to \hat{a}$ is a 
$C_{0}(X)$-isomorphism from $A$ onto the $C_{0}(X)$-algebra $C_{0}(X,A^{\sharp})$ of continuous sections of $A^{\sharp}$ that vanish at infinity. 
\end{theorem}

If $A, B$ are $C_{0}(X)$-algebras and $T:A\to B$ is a 
$C_{0}(X)$-morphism, then $T^{\sharp}:A^{\sharp}\to B^{\sharp}$ is continuous, where 
$T^{\sharp}a_{x}=T_{x}(a_{x})=(Ta)_{x}$.  (In fact $(T^{\sharp})^{-1}(U(Ta,\eps))\supset U(a,\eps)$.)  

Next, if $\mfB$ is a $C_{0}(X)$-subalgebra of a $C_{0}(X)$-algebra $B$ then for any $x$, $I_{x}B\cap \mfB=I_{x}\mfB$ so that we can identify $\mfB^{\sharp}$ with a subbundle of $B^{\sharp}$, and the topology on $\mfB^{\sharp}$ is the relative topology.

We now recall how the (maximal) tensor product $A\otimes_{C_{0}(X)} B$ of two $C_{0}(X)$-algebras is defined.  For more information, see \cite{Blanchard92,Blanchard95,LeGall,EchterWill}.  One natural way to do this is to take
$A\otimes_{C_{0}(X)} B$
to be the maximal $C_{0}(X)$-balanced tensor product: so $A\otimes_{C_{0}(X)} B=(A\otimes_{max} B)/I$, where $I$ is the closed ideal generated by differences of the form $(af\otimes b - a\otimes fb)$ ($a\in A, b\in B, f\in C_{0}(X)$). 
The $C_{0}(X)$-action on $A\otimes_{C_{0}(X)} B$ is determined by: 
$f(a\otimes b)=fa\otimes b=a\otimes fb$ for $f\in C_{0}(X)$.  (Alternatively, one regards $A\otimes_{max} B$ as a 
$C_{0}(X\x X)$-algebra and ``restricts to the diagonal'': $A\otimes_{C_{0}(X)} B=
(A\otimes_{max} B)/C_{\De}(A\otimes_{max} B)$  where $C_{\De}=\{g\in C_{0}(X\x X): g(x,x)=0 \mbox{ for all }
x\in X\}$.)  Next, $(A\otimes_{C_{0}(X)} B)_{x}=A_{x}\otimes_{max} B_{x}$.  If $D$ is an ordinary $\cs$, then 
$D\otimes_{max} B$ is a $C_{0}(X)$-algebra in the natural way: $\theta(f)(d\otimes b)=d\otimes fb$.  (Alternatively, one can identify the $C_{0}(X)$-algebra $D\otimes_{max} B$ with $(D\otimes_{max} C_{0}(X))\otimes_{C_{0}(X)} B$.)  
An important case of this is when 
$D=C_{0}(Z)$ ($Z$ a locally compact Hausdorff space): then $C_{0}(Z,B)=C_{0}(Z)\otimes B$ is a $C_{0}(X)$-algebra, and 
$C_{0}(Z,B)_{x}=C_{0}(Z)\otimes B_{x}=C_{0}(Z,B_{x})$.  It is easily checked that $(g\otimes b)_{x}=g\otimes b_{x}$, and it follows that for $F\in C_{0}(Z,B)$, $F_{x}(z)=F(z)_{x}\in B_{x}$.    

Now let $B=C_{0}(Y)$ ($Y$ a locally compact Hausdorff space) with the $C_{0}(X)$-action on $B$ given by a continuous map 
$q:Y\to X$: here $(fF)(y)=f(q(y))F(y)$ where $F\in C_{0}(Y), f\in C_{0}(X)$.  In this case, one writes
$q^{*}A=A\otimes_{C_{0}(X)} C_{0}(Y)$.  It is sometimes helpful to incorporate explicit mention of the map $q$ in this tensor product by writing $A\otimes_{C_{0}(X),q} C_{0}(Y)$ in place of 
$A\otimes_{C_{0}(X)} C_{0}(Y)$.  $q^{*}A$ is actually also a $C_{0}(Y)$-algebra in the obvious way: $(a\otimes F)F'=a\otimes FF'$ for $F,F'\in C_{0}(Y)$, and for each $y\in Y$, we have 
$(q^{*}A)_{y}=A_{q(y)}$.  (The canonical map from $(q^{*}A)_{y}$ to $A_{q(y)}$ comes from sending 
$(a\otimes F)_{y}$ to $a_{q(y)}F(y)$.)  
Now let $Y\x_{q} A^{\sharp}=\{(y,a_{q(y)}): y\in Y, a\in A\}$ with the relative topology inherited from 
$Y\x A^{\sharp}$.  From the above, $Y\x_{q} A^{\sharp}$ is identified {\em as a set} with $(q^{*}A)^{\sharp}$.  We now show that the spaces are homeomorphic when $q$ is open.

\begin{proposition}   \label{prop:p*A}
If $q$ is also open, then the identity map $i:Y\x_{q} A^{\sharp}\to (q^{*}A)^{\sharp}$ is a homeomorphism.
\end{proposition}
\begin{proof}
A base for the topology of $(q^{*}A)^{\sharp}$ is given by sets of the form $W(a\otimes F,\eps)$ where $W$ is a relatively compact open subset of $Y$, $F\in C_{0}(Y)$ is $1$ on $W$ and $a\in A$.  Then 
$W(a\otimes F,\eps)=W\x_{q} q(W)(a,\eps)$ and the latter sets form a base for the topology of $Y\x_{q} A^{\sharp}$.  
\end{proof}

We now recall what is meant by a $\Ga$-algebra $A$ (\cite{LeGall,Popescu}).  Form the balanced tensor 
products $s^{*}A=A\otimes _{C_{0}(X),s} C_{0}(\Ga)$ and $r^{*}A=A\otimes _{C_{0}(X),r} C_{0}(\Ga)$.  From 
Theorem~\ref{th:asharp}, $r^{*}A=C_{0}(\Ga,(r^{*}A)^{\sharp})$.
Then $A$ is called a 
{\em $\Ga$-algebra} if there is given a $C_{0}(\Ga)$-isomorphism $\al:s^{*}A\to r^{*}A$ such that the induced isomorphisms 
$\al_{\ga}:(s^{*}A)_{\ga}=A_{s(\ga)} \to (r^{*}A)_{\ga}=A_{r(\ga)}$ satisfy the groupoid multiplication properties:
$\al_{\ga\ga'}=\al_{\ga}\al_{\ga'}$ whenever $r(\ga')=s(\ga)$
and $\al_{\ga^{-1}}=(\al_{\ga})^{-1}$.  Obviously, for each $x\in X$, $\al_{x}$ is the identity map on $A_{x}$.

As an example, suppose that $\Ga$ is a locally compact group $G$.  Then $s^{*}A=r^{*}A=A\otimes C_{0}(G)=C_{0}(G,A)$.  
For $F\in C_{0}(G,A)$, we have $F_{g}=F(g)$ and by Proposition~\ref{prop:p*A}, 
$(C_{0}(G,A))^{\sharp}=G\x A$.  If $\al:C_{0}(G,A)\to C_{0}(G,A)$ gives a $G$-action on $C_{0}(G,A)$,
then since $(C_{0}(G,A))_{g}=A$, we get isomorphisms
$\al_{g}:A\to A$.  We are given that $\al_{gh}=\al_{g}\circ \al_{h}$ for all $g,h\in G$.  Last since the map
$g\to (\al(a\otimes k))_{g}=\al_{g}(a)k(g)$ 
belongs to $C_{0}(G,A)$, it follows that for each $a\in A$, the map 
$g\to \al_{g}(a)$ is norm continuous. So $A$ is a $G$-algebra in the usual sense.
(The converse is left to the reader.)  We now show that the groupoid version of the preceding holds; a 
$\Ga$-$\cs$ $A$ can then be viewed in terms of a $\Ga$-action on $A^{\sharp}$. 
Note that the corollary characterizes $\Ga$-action in terms exactly analogous to that of a group action.

\begin{theorem}   \label{th:gacs}
$A$ is a $\Ga$-algebra if and only if $A^{\sharp}$ is a continuous $\Ga$-space of $\css$ $A_{x}$.  
\end{theorem}
\begin{proof}
Suppose that $A$ is a $\Ga$-algebra.  So we are given a $C_{0}(\Ga)$-isomorphism $\al:s^{*}A\to r^{*}A$ with 
$\ga\to \al_{\ga}$ a homomorphism into $Iso(A^{\sharp})$.  Let $\bt=\al^{\sharp}$.  Then the continuity of the map 
$(\ga,z)\to \al_{\ga}(z)$ follows, using 
Proposition~\ref{prop:p*A}, by composing the following continuous maps:\\
\[  \Ga\x_{s} A^{\sharp}\overset{i}{\to} (s^{*}A)^{\sharp}\overset{\bt}{\to} 
(r^{*}A)^{\sharp}\overset{i^{-1}}\to \Ga\x_{r} A^{\sharp}\overset{p_{2}}\to A^{\sharp}  \]
where $p_{2}$ is the projection onto the second coordinate.  So  $A^{\sharp}$ is a continuous $\Ga$-space of $\css$.

Conversely, suppose that $A^{\sharp}$ is a continuous $\Ga$-space of $\css$.  Define 
$\bt:\Ga\x_{s} A^{\sharp}\to \Ga\x_{r} A^{\sharp}$ by: $\bt(\ga,z)=(\ga,\al_{\ga}(z))$.  By assumption, $\bt$ is continuous.  The map $\bt^{-1}$ is also continuous since it equals 
$(inv\otimes 1)\circ \bt\circ (inv\otimes 1)$, where $inv(\ga)=\ga^{-1}$.  Then the map 
$F\to (i\bt i^{-1})\circ F$ is a $C_{0}(\Ga)$-homomorphism from $C_{0}(\Ga,(s^{*}A)^{\sharp})$ into 
$C_{0}(\Ga,(r^{*}A)^{\sharp})$ which is an isomorphism since it's inverse is the corresponding expression involving  
$\bt^{-1}$.  By Theorem~\ref{th:asharp}, this isomorphism determines a $C_{0}(\Ga)$-isomorphism $\al:s^{*}A\to r^{*}A$ and $A$ is a $\Ga$-algebra.
\end{proof}
\begin{corollary}   \label{cor:gacs}
$A$ is a $\Ga$-algebra if and only if there is given a groupoid homomorphism $\ga\to \al_{\ga}$ from $\Ga$ into 
$Iso(A^{\sharp})$ such that for each $a\in A$, the map $\ga\to \al_{\ga}(a_{s(\ga)})$ is continuous.
\end{corollary}
\begin{proof}
Suppose that we are given a groupoid homomorphism $\ga\to \al_{\ga}$ from $\Ga$ into 
$Iso(A^{\sharp})$ such that for each $a\in A$, the map $\ga\to \al_{\ga}(a_{s(\ga)})$ is continuous.  Let 
$\{(\ga_{\de},z_{\de})\}$ be a net in $\Ga \x_{s} A^{\sharp}$ converging to some $(\ga_{0},z_{0})$.  We show that 
$\al_{\ga_{\de}}(z_{\de})\to \al_{\ga_{0}}(z_{0})$.  Let $a, c\in A$ be such that $z_{0}=a_{s(\ga_{0})}$ and
$\al_{\ga_{0}}(z_{0})=c_{r(\ga_{0})}$.  Let $V(c,\eps)$ be a neighborhood of $c_{r(\ga_{0})}$ in $A^{\sharp}$.
By continuity of the map $\ga\to \al_{\ga}(a_{s(\ga)})$, there exists an open neighborhood $Z$ of $\ga_{0}$ in $\Ga$
and a $\de_{1}$ such that for all $\de\geq \de_{1}$, $\al_{\ga}(a_{s(\ga)})\in V(c,\eps/2)$ for all $\ga\in Z$.  Since 
$z_{\de}\to z_{0}$, we can also arrange that $z_{\de}\in s(Z)(a,\eps/2)$ for all $\de\geq \de_{1}$.  So for all 
$\de\geq \de_{1}$, $\norm{z_{\de} - a_{s(\ga_{\de})}}<\eps/2$, giving 
$\norm{\al_{\ga_{\de}}(z_{\de}) - \al_{\ga_{\de}}(a_{s(\ga_{\de})})}<\eps/2$.   Since 
$\al_{\ga_{\de}}(a_{s(\ga_{\de})})\to c_{r(\ga_{0})}$, we can also suppose that for all $\de\geq \de_{1}$,  
$\al_{\ga_{\de}}(a_{s(\ga_{\de})})\in V(c,\eps/2)$, i.e. 
$\norm{\al_{\ga_{\de}}(a_{s(\ga_{\de})}) - c_{r(\ga_{\de})}}<\eps/2$.   By the triangular inequality, 
$\al_{\ga_{\de}}(z_{\de})\in V(c,\eps)$ ($\de\geq \de_{1}$), and $\al_{\ga_{\de}}(z_{\de})\to \al_{\ga_{0}}(z_{0})$.  By Theorem~\ref{th:gacs}, $A$ is a $\Ga$-algebra.  The converse also follows from Theorem~\ref{th:gacs}.
\end{proof}

Now suppose that $A$ is a $\Ga$-algebra and $J$ is a closed ideal of $A$ that is a $\Ga$-subalgebra of $A$ in the natural way, i.e. for each $\ga \in \Ga$, $j\in J$, we 
have $\al_{\ga}(j_{s(\ga)})\in J_{r(\ga)}$.  Using the continuity of the canonical map from $A^{\sharp}$ to 
$(A/J)^{\sharp}$ and Corollary~\ref{cor:gacs}, it is easy to prove that the $C_{0}(X)$-algebra $A/J$ is also a 
$\Ga$-algebra in the natural way, and we have a short exact sequence of $\Ga$-algebras:
\begin{equation}  \label{eq:ses}
0\to J\to A\to A/J\to 0.
\end{equation}

Next suppose that $A$ is a $C_{0}(X)$-algebra, and that the $A_{x}$'s form a $\Ga$-space of $\css$.  So we can say that 
$A$ has an {\em algebraic} $\Ga$-action (with no continuity condition on the maps 
$\ga\to \al_{\ga}(a_{s(\ga)})$).  We wish to define a $C^{*}$-subalgebra $A^{cont}$ of $A$ on which the $\Ga$-action
{\em is} continuous.  For this result, we require that $\Ga$ have {\em local $r-G$-sets}
(cf. \cite[p.10]{rg}, \cite[p.44]{Paterson}).  This means that for each $\ga_{0}\in \Ga$, there exists an open neighborhood $U$ of 
$r(\ga_{0})$ in $X$ and a subset $W$ of $\Ga$ containing $\ga_{0}$ such that $r_{W}=r_{\mid W}$ is a homeomorphism from $W$ onto 
$U$.  Most locally compact groupoids that arise in practice have local $r-G$-sets (e.g. Lie groupoids, r-discrete groupoids and transformation group groupoids).  

\begin{proposition}   \label{prop:gacs}
Let $\Ga$ have local $r-G$-sets, $A$ have an algebraic $\Ga$-action and define
\[  A^{cont}=\{a\in A: \mbox{ the map } \ga\to \al_{\ga}(a_{s(\ga)}) \mbox{ is continuous}\}.  \]
Then $A^{cont}$ is a $\Ga$-subalgebra of $A$.
\end{proposition}
\begin{proof}
It is obvious from the definition of a $C^{*}$-bundle that $A^{cont}$ is a $^{*}$-subalgebra of $A$.  Let $a_{n}\to a$ in $A$ with $a_{n}\in A^{cont}$ for all $n$.  Then $\al_{\ga}((a_{n})_{s(\ga)})\to \al_{\ga}(a_{s(\ga)})$ 
uniformly in $\ga$.  Adapting the proof of the elementary result that a uniform limit of continuous functions is continuous - one uses also the upper semicontinuity of the norm on $A^{\sharp}$ - it follows that the map
$\ga\to \al_{\ga}(a_{s(\ga)})$ is continuous, i.e. $a\in A^{cont}$.  So $A^{cont}$ is a $C^{*}$-subalgebra of $A$.
Next, if $f\in C_{0}(X)$, then $\al_{\ga}((fa)_{s(\ga))})=f(s(\ga))\al_{\ga}(a_{s(\ga))})$, and the map
$\ga\to \al_{\ga}((fa)_{s(\ga)})$ is continuous.  So $C_{0}(X)A^{cont}\subset A^{cont}$.  Also, if $a\in A^{cont}$,
then $a=f'a'$ for some $f'\in C_{0}(X), a'\in A$, and so $a=\lim e_{n}(f'a')$ where $\{e_{n}\}$ is a bounded approximate identity for $C_{0}(X)$.  So $C_{0}(X)A^{cont}=A^{cont}$, and $A^{cont}$ is a $C_{0}(X)$-algebra.  Last, we have to show that if $a\in A^{cont}$ and $\ga_{0}\in \Ga$, then $\al_{\ga_{0}}(a_{s(\ga_{0})})\in (A^{cont})_{r(\ga_{0})}$.  Let $W, U$ be as above so that $r_{W}:W\to U$ is a homeomorphism.  Let $f\in C_{c}(U)$ be such that $f(r(\ga_{0}))=1$.  Then the section $g$ of $A^{\sharp}$ given by: $g(u)=f(u)\al_{r_{W}^{-1}(u)}(a_{s(r_{W}^{-1}(u))})$ belongs to 
$C_{c}(X,A^{\sharp})$.  By Theorem~\ref{th:asharp}, there exists $b\in A$ such that $b_{u}=g(u)$ for all $u\in X$.  
Since $\al_{\ga}(b_{s(\ga)})=
f(s(\ga))\al_{\ga\cdot r_{W}^{-1}(s(\ga))}(a_{s(\ga\cdot r_{W}^{-1}(s(\ga)))})$ and $a\in A^{cont}$, we see that 
$b\in A^{cont}$.
\end{proof}  

Now let $A, B$ be $\Ga$-algebras.  The tensor product $A\otimes_{C_{0}(X)} B$ is a $\Ga$-algebra (\cite[3.1.2]{LeGall}) in a natural way.  Indeed, using the associativity of the balanced tensor product (\cite[p.90]{Blanchard92}) and the equality
$C_{0}(\Ga)\otimes_{C_{0}(\Ga)} C_{0}(\Ga)=C_{0}(\Ga)$, we obtain 
$q^{*}(A\otimes_{C_{0}(X)} B)=q^{*}A\otimes_{C_{0}(\Ga)} q^{*}B$ ($q=s,r$).  The $\Ga$-action on $A\otimes_{C_{0}(X)} B$ 
is then given by $\al\otimes \bt$, where $\al, \bt$ are the $\Ga$-actions on $A,B$.   
Further $(\al\otimes \bt)_{\ga}=\al_{\ga}\otimes \bt_{\ga}$ (recalling that 
$(A\otimes_{C_{0}(X)} B)_{x}=A_{x}\otimes_{max} B_{x}$). 
Also, if $A$ is just a $\cs$ and $B$ is a $\Ga$-algebra, then the $C_{0}(X)$-algebra $A\otimes_{max}B$ 
is a $\Ga$-algebra: we identify $q^{*}(A\otimes_{max} B)$ with $A\otimes_{max} q^{*}B$ with $q=s,r$, and the $\Ga$-action is given by $I\otimes \bt$.  (Alternatively, one can reduce this to the earlier case by using 
$q^{*}((A\otimes C_{0}(X))\otimes_{C_{0}(X)} B)$.)  In particular, if $B$ is a $\Ga$-algebra, then the $C_{0}(X)$-algebra
$C_{0}(T,B)$ is also a $\Ga$-algebra, and the action is given by:
\begin{equation}  \label{eq:gactionc0}
\al_{\ga}(F_{s(\ga)})(t)=\al_{\ga}(F(t)_{s(\ga)}).
\end{equation}  

A $\Ga$-homomorphism (\cite[Definition 3.1.2]{LeGall}) from $A$ to $B$ is a $C_{0}(X)$-homomorphism 
$\phi:A\to B$ such that for all $\ga\in \Ga$, 
\begin{equation}  \label{eq:Gahomo}
\phi_{r(\ga)}\al_{\ga}=\bt_{\ga}\phi_{s(\ga)}
\end{equation}
where $\al, \bt$ denote respectively the actions of $\Ga$ on $A$ and $B$.  It is simple to check that with $\Ga$-homomorphisms as morphisms, the class of $\Ga$-algebras forms a category.


\section{Continuity and exactness}

Next, we need the notion of a {\em crossed product} of $\Ga$ by a $\Ga$-algebra $A$ (\cite{Ren87}, \cite{LeGall}).  We will need to use the profound disintegration theorem of J. Renault of \cite{Ren87}.
Renault develops a groupoid version of the theory of twisted covariance algebras for locally compact groups, and working in a very general context, constructs a $\cs$
$C^{*}(\Ga,\Si,A,\la)$ where $A$ is a $\Ga$-algebra and $S$ is a bundle of abelian groups over $X$ with $\Ga$ acting on the fibers and $\Si$ is a groupoid given by an exact sequence of groupoids.  (Also used in the construction is a homomorphism 
$\chi$ on $S$ that we don't need to go into for our present purposes.)  For the special case of the groupoid crossed product, we take $S=X=\Ga^{0}$ and 
$\Si=\Ga$.  In that context, we put on the algebra $r^{*}_{c}A=C_{c}(\Ga,(r^{*}A)^{\sharp})\cong C_{c}(\Ga)r^{*}A$, a product and involution given by:
\begin{equation}  \label{eq:conv}
 F_{1}*F_{2}(\ga)=\int_{\Ga^{r(\ga)}}F_{1}(\ga')\al_{\ga'}(F_{2}(\ga'^{-1}\ga))\,d\la^{r(\ga)}(\ga')
\hspace{.2in} (F_{1})^{*}(\ga)=\al_{\ga}(F_{1}(\ga^{-1})^{*}).  
\end{equation}
(The proof that $F_{1}*F_{2}\in r_{c}^{*}A$ is given by P.-Y. Le Gall in \cite[Proposition 7.1.1]{LeGall}.)  Next,
$r^{*}_{c}A$ is a normed $^{*}$-algebra with isometric involution under the $I$-norm $\norm{.}_{I}$, where 
\[   \norm{F}_{I}=\max\{\norm{F}_{r}, \norm{F^{*}}_{r}\}              \]
and 
\[   \norm{F}_{r}=\sup_{x\in X} \int_{\Ga^{x}}\norm{F(\ga)}\,d\la^{x}(\ga).       \]
The enveloping $\cs$ of $(r^{*}_{c}A,\norm{.}_{I})$ is then defined to be the crossed product 
$C^{*}(\Ga,A)$.

A very simple example of this is provided by the case where $A=C_{0}(X)$ with the usual action of 
$\Ga$ on $X$: $\al_{\ga}(s(\ga))=r(\ga)$.  In that case, as is easily checked, $A^{\sharp}=X\x \C$,
$\al_{\ga}:\C_{s(\ga)}\to \C_{r(\ga)}$ is the identity map, $s^{*}A=r^{*}A=C_{0}(\Ga)$, 
and $\al: C_{0}(\Ga)\to C_{0}(\Ga)$ the identity map.  Of course,
$(r^{*}A)^{\sharp}$ is just $\Ga\x \C$, and $r^{*}_{c}A=C_{c}(\Ga)$.  Then
$C^{*}(\Ga,A)$ is just the $\cs$ $C^{*}(\Ga)$ of the groupoid (\cite{rg}).

We now turn to Renault's disintegration theorem for representations of $C^{*}(\Ga,A)$ - for a detailed exposition for the case $C^{*}(\Ga)$ see \cite{MuhlyTCU}.  The theory uses the fundamental papers of Ramsay
(\cite{Ramvirt,Ramtop}).  We first formulate \cite[Lemme 4.5]{Ren87} in $C_{0}(X)$-algebra terms.  Let $A$ be a $C_{0}(X)$-algebra, $\mathfrak{H}=\{H_{x}\}_{x\in X}$ a Hilbert bundle and $\mu$ a probability measure on $X$.  
Let $H=L^{2}(X,\mu,\mathfrak{H})$.  We will say that a non-degenerate representation $\pi:A\to B(H)$ is a 
{\em $C_{0}(X)$-representation} (for $(X,\mu,\mathfrak{H})$)
if $\pi$ commutes with the $C_{0}(X)$-actions on $A$ and $B(H)$, i.e. for all $f\in C_{0}(X)$ and all
$a\in A$, $T_{f}\pi(a)=\pi(a)T_{f}=\pi(fa)$, where $T_{f}$ is the multiplication operator on 
$L^{2}(H)$ associated with $f$.  By taking strong operator limits, we get that every $a$ commutes with every $T_{f}$ for 
$f\in L^{\infty}(X,\mu)$, i.e. with every diagonalizable operator.  So (\cite[II, 2, 5, Corollary]{D1}) every $\pi(a)$ is decomposable, and from \cite[Lemma 8.3.1]{D2}, $\pi$ is a direct integral of representations $\pi_{x}$ of $A$.  Further,
for each $f\in C_{0}(X)$, $\pi_{x}(fa)=f(x)\pi_{x}(a)$ so that $\pi_{x}$ is a representation of $A_{x}$ on $H_{x}$.  The $\pi_{x}$'s are non-degenerate a. e. by \cite[8.1.5]{D2}.  We now discuss what is meant by a covariant representation of
$(\Ga,A)$.

Let $\mu$ be quasi-invariant on $X$, $\nu=\int_{X}\la^{x}\,d\mu(x)$ (\cite[pp.22-23]{rg}, \cite[Ch. 3]{Paterson}): quasi-invariance means that $\nu\sim \nu^{-1}$.  Let $U$ be a Borel subset of $X$ which is $\mu$-conull.   Then 
$\Ga_{\mid U}=r^{-1}(U)\cap s^{-1}(U)$ is a Borel groupoid which is $\nu$-conull in $\Ga$.  Then 
$\Ga_{\mid U}$ equipped with the restrictions of $\mu, \nu$ to $U, \Ga_{\mid U}$ is a measured groupoid, called the 
{\em inessential contraction} of $\Ga$ to $U$.   Next we are given a Hilbert bundle $\mathfrak{H}=\{H_{x}\}_{x\in X}$; 
$Iso(X*\mathfrak{H})$ is the Borel groupoid of unitaries $U_{y,x}:H_{x}\to H_{y}$ as $x,y$ range over $X$. (See \cite[Chapter 3]{MuhlyTCU}.)  A {\em covariant representation} (or {\em a representation of the dynamical system $(\Ga,\Ga,A)$}
in the terminology of (\cite[p.79]{Ren87})) $(L,\pi)$ of the pair $(\Ga,A)$ consists of:
\bi
\item[(i)] a Borel homomorphism $L:\Ga_{\mid U}\to Iso(X*\mathfrak{H})_{\mid U}$,
\item[(ii)]  a (non-degenerate) $C_{0}(X)$-representation $\pi=\int^{\oplus}\pi_{x}\,d\mu(x)$ 
of $A$ on $H=L^{2}(X,\mu,\mfH)$: for each $a\in A$, 
\[  \pi(a)=\int^{\oplus} \pi_{x}(a_{x})\,d\mu(x), \]
\item[(iii)] for all $\ga\in \Ga_{\mid U}$ and $a\in A$, we have
\begin{equation}   \label{eq:covar}
L_{\ga}\pi_{s(\ga)}(a_{s(\ga)})L_{\ga^{-1}}=\pi_{r(\ga)}(\al_{\ga}(a_{s(\ga)})).
\end{equation}
\ei
A. Ramsay (\cite{Ramvirt,Ramtop,MuhlyTCU}) showed, at least in the case $A=C_{0}(X)$ above we can actually take $U=X$.  However, because of the conullity of $U$, we can effectively regard the pair $(L,\pi)$ as defined on $\Ga$ rather than 
$\Ga_{\mid U}$ and will usually leave the $U$ implicit. 

Every covariant representation $(L,\pi)$ of $\Ga$ integrates up to give a representation $\Phi$ of $C^{*}(\Ga,A)$.  Indeed, from \cite[p.80]{Ren87}, for $F\in r^{*}_{c}A$ and $\xi,\eta\in \mathfrak{H}$, 
\begin{equation}   \label{eq:integform}
\lan \Phi(F)\xi,\eta\ran    
=\int\lan\pi_{r(\ga)}(F_{\ga})L_{\ga}\xi_{s(\ga)},\eta_{r(\ga)}\ran\,d\nu_{0}(\ga)   
\end{equation}
where, as usual (\cite[p.52]{rg}, \cite[3.1]{Paterson}) $d\nu_{0}=D^{-1/2}\,d\nu$ with $D=d\nu/d\nu^{-1}$.  

Conversely, every representation $\Psi$ of $C^{*}(\Ga,A)$ on a Hilbert space $K$ 
is equivalent to such a $\Phi$.  Indeed, from \cite[p.88]{Ren87}, elements $\phi, h$ of 
the algebras $C_{c}(\Ga)$, $C_{c}(X,A^{\sharp})$ act as left multipliers on $r^{*}_{c}A$ where:
\begin{equation}  \label{eq:mult}
\phi*F(\ga)=\int_{\Ga^{r(\ga)}}\phi(\ga')\al_{\ga'}(F(\ga'^{-1}\ga))\,\la^{r(\ga)}(\ga'), \hspace{.2in}
(hF)(\ga)=(h\circ r)(\ga)F(\ga).
\end{equation}
Renault shows that there are representations $L', \pi'$ of $C_{c}(\Ga), C_{c}(X, A^{\sharp})=C_{c}(X)A$ on 
$K$ and determined by:
\begin{equation}   \label{eq:mult'}
\Psi(\phi*F)=L'(\phi)\Psi(F), \hspace{.2in}  \Psi(hF)=\pi'(h)\Psi(F).               
\end{equation}
Renault first studies the representation $L'$ of $C_{c}(\Ga)$ and obtains a quasi-invariant measure $\mu$ on $X$ and a measurable Hilbert bundle $\mfH =\{H_{x}\}$ over $X$ such that $K$ can be identified with $L^{2}(X,\mu;\mfH)$.  Then $L'$ is disintegrated into a representation $L$ of the groupoid $\Ga$, and the
$\pi'(a)$'s are decomposable on $\mfH$: $\pi'(a)=\{\pi_{x}(a_{x})\}$ (with $\pi$ in place of $\pi'$).  He then shows that $L,\pi$ can be taken to be such that the pair $(L,\pi)$ is a covariant pair whose integrated form is equivalent to $\Psi$.

We now discuss exactness and continuity for functors.  So let
${\bf F}$ be a functor from the category of $\Ga$-algebras with $\Ga$-homomorphisms as morphisms into the category of ordinary 
$\css$ with $^{*}$-homomorphisms as morphisms.  Following \cite[p.19, ff.]{GHT}, we say that ${\bf F}$ is 
{\em exact} if for every short exact sequence of $\Ga$-algebras
\[  0\to J\to A\to A/J\to 0   \]
the induced sequence of ordinary $\css$
\[   0\to {\bf F}(J)\to {\bf F}(A)\to {\bf F}(A/J)\to 0          \]
is exact.  Now let $I$ be a closed interval $[a,b]$, $B$ be a $\Ga$-algebra and $IB$ be the $\Ga$-algebra 
$C(I)\otimes B=C(I,B)$.  For each $k\in {\bf F}(IB)$, we can associate a function 
$\hat{k}:I\to {\bf F}(B)$ by setting 
$\hat{k}(t_{0})={\bf F}(ev_{t_{0}})(k)$ where $ev_{t_{0}}:IB\to B$ is evaluation at $t_{0}$: 
$ev_{t_{0}}(g)=g(t_{0})$.  (Note that $ev_{t_{0}}$ is a $\Ga$-homomorphism.)  The functor ${\bf F}$ is said to be 
{\em continuous} if every $\hat{k}$ is continuous.  The map $k\to \hat{k}$ then gives a homomorphism from ${\bf F}(IB)$ into $I{\bf F}(B)$.   Later, we will need to replace the finite interval $I$ in $C(I,B)$ by the infinite interval $T$.  We cannot replace $C(I,B)$ by $C_{b}(T,B)$ since the latter does not admit a $\Ga$-action in any natural way.  However, the theory can be made to work, as we will see later (Proposition~\ref{prop:Tcont}) by replacing 
$C_{b}(T,B)$ by a $C_{0}(X)$-subalgebra $\mfB$ with a covering $\Ga$-action. 

An exact, continuous functor ${\bf F}$ will now be constructed from the category of $\Ga$-algebras with 
$\Ga$-homomorphisms as morphisms into the category of ordinary $\css$ with $^{*}$-homomorphisms as morphisms.  For a 
$\Ga$-algebra $A$, define 
${\bf F}(A)=C^{*}(\Ga,A)$.  We need to specify what ${\bf F}$ does to morphisms.   Let $B$ also be a $\Ga$-algebra and 
$\phi:A\to B$ be a $\Ga$-homomorphism.  Define, for each
$F\in r^{*}_{c}A$, a section $\Tilde{\phi}(F):\Ga\to (r^{*}B)^{\sharp}$ by:
\begin{equation}  \label{eq:phitild}
\Tilde{\phi}(F)(\ga)=\phi_{r(\ga)}(F_{\ga}).  
\end{equation}
We note that $\Tilde{\phi}(F)$ is just the same as 
$(r^{*}\phi)(F)\in r_{c}^{*}B=C_{c}(\Ga,(r^{*}B)^{\sharp})$.  
Using (\ref{eq:conv}) and (\ref{eq:Gahomo}), we obtain that for $F_{1},F_{2}\in r^{*}_{c}A$ and each $\ga\in \Ga$,
\begin{equation} \label{eq:star}
\Tilde{\phi}(F_{1}*F_{2})(\ga)=\phi_{r(\ga)}((F_{1}*F_{2})(\ga))=(\Tilde{\phi}(F_{1})*\Tilde{\phi}(F_{2}))(\ga)   
\end{equation}
and $\Tilde{\phi}((F_{1})^{*})=(\Tilde{\phi}(F_{1}))^{*}$.  
So $\Tilde{\phi}$ is a $^{*}$-homomorphism from $r^{*}_{c}A$ to 
$r^{*}_{c}B$.  It is continuous for the respective $C^{*}$-norms since $\norm{\Tilde{\phi}(F)}_{I}\leq \norm{F}_{I}$ so that $\pi\circ \Tilde{\phi}$ is a  representation of $C^{*}(\Ga,A)$ whenever $\pi$ is a representation of 
$C^{*}(\Ga,B)$.  We set ${\bf F}(\phi)=\Tilde{\phi}$.  It is easy to check that ${\bf F}$ is a functor.  In the following, $A\otimes B$ means $A\otimes_{max} B$.

\begin{theorem}  \label{th:exact}
The functor ${\bf F}$ is continuous and exact.
\end{theorem}
\begin{proof}
(a) We first show the continuity of $F$ (cf. \cite[Lemma 4.11]{GHT} where the locally compact group case is sketched).   Let $A$ be an ordinary $\cs$ and $B$ be a $\Ga$-$\cs$ and recall that 
$A\otimes B$ is a $\Ga$-algebra.  We show that 
\begin{equation}  \label{eq:equal}
A\otimes C^{*}(\Ga,B)\cong C^{*}(\Ga,A\otimes B).  
\end{equation}
There is a natural $C_{0}(\Ga)$-isomorphism $\Phi$ from $A\otimes r^{*}B$ onto 
$r^{*}(A\otimes B)$ determined by:
$\Phi(a\otimes (b\otimes_{C_{0}(X)} F))=(a\otimes b)\otimes_{C_{0}(X)} F$ (e.g. \cite[Corollary 1.3]{Popescu}).  
The map $\Phi$ restricts to an isomorphism, also denoted $\Phi$, from $A\otimes_{alg} r_{c}^{*}B$ onto a subalgebra of 
$r_{c}^{*}(A\otimes B)$.  $\Phi$ is also an isomorphism when $A\otimes_{alg} r_{c}^{*}B, r_{c}^{*}(A\otimes B)$ are given the convolution product and involution of (\ref{eq:conv}).   Give 
$A\otimes_{alg} r_{c}^{*}B$, $r_{c}^{*}(A\otimes B)$ the norms that they inherit as (dense) subalgebras of 
$A\otimes C^{*}(\Ga,B)$, $C^{*}(\Ga,A\otimes B)$.  We note that $\Phi(A\otimes r_{c}^{*}B)$ is $\norm{.}_{I}$-dense in 
$r_{c}^{*}(A\otimes B)$ and so also dense in $C^{*}(\Ga,A\otimes B)$.  We show that $\Phi$ is isometric.

A representation $\pi'$ of $r_{c}^{*}(A\otimes B)$ is determined by a covariant pair $(L,\pi)$ where $\pi$ is a representation of $A\otimes B$ on some $H=L^{2}(X,\mu,\mfH)$.  Then (\cite[Theorem 6.3.5]{Murphy}) there exist non-degenerate, commuting representations $\pi_{1}, \pi_{2}$ of $A, B$ on $H$ such that 
$\pi(a\otimes b)=\pi_{1}(a)\pi_{2}(b)$.  Further, using bounded approximate identities in $A, B$, $\pi_{2}$ is a 
$C_{0}(X)$-representation and $\pi_{1}$ commutes with the $C_{0}(X)$-multiplication operators $T_{f}$ on $H$.  
Disintegrating, we get 
\[   \pi_{1}(a)_{x}(\pi_{2})_{x}(b_{x})=\pi_{x}(a\otimes b_{x})=(\pi_{2})_{x}(b_{x})\pi_{1}(a)_{x}   \]
almost everywhere, and $L_{\ga}\pi_{1}(a)_{s(\ga)}(\pi_{2})_{s(\ga)}(b_{s(\ga)})L_{\ga^{-1}}=
\pi_{1}(a)_{r(\ga)}(\pi_{2})_{r(\ga)}(b_{r(\ga)})$.  It follows that $(L,\pi_{2})$ is covariant for $B$ and
$L_{\ga}\pi_{1}(a)_{s(\ga)}L_{\ga^{-1}}=\pi_{1}(a)_{r(\ga)}$ a.e..  Let $\Phi_{2}$ be the integrated form of 
$(L,\pi_{2})$.  Then from (\ref{eq:integform}) and the above, the representations $\pi_{1}, \Phi_{2}$ commute, 
and so the $C^{*}$-semi-norm that they induce on $A\otimes_{alg} C^{*}(\Ga,B)$ is $\leq$ the
maximum tensor product norm.  Since $\pi'(\Phi(w))=\\
(\pi_{1}\otimes \Phi_{2})(w)$ ($w\in A\otimes_{alg} r_{c}^{*}B$), it follows that $\norm{\Phi(w)}\leq \norm{w}$.

On the other hand, each representation $\pi$ of $A\otimes C^{*}(\Ga,B)$ is determined by a pair of commuting representations $\pi_{1}, \pi_{2}$ of $A, C^{*}(\Ga,B)$ on some $H$.  Then $\pi_{2}$ disintegrates into a 
covariant pair $(\pi_{2}', L)$ and we can identify $H=L^{2}(X,\mu,\mfH)$.  Using (\ref{eq:mult'}), $\pi_{1}$ and 
$\pi_{2}'$ commute and $L'$ and 
$\pi_{1}$ commute.  Also, $\pi_{1}$ commutes with the $T_{f}$'s ($f\in C_{0}(X)$).  From the proof of the disintegration theorem, $L_{\ga}\pi_{1}(a)_{s(\ga)}L_{\ga^{-1}}=\pi_{1}(a)_{r(\ga)}$ almost everywhere.  Then
$(\pi_{1}\otimes \pi_{2}', L)$ is a covariant representation for $A\otimes B$ and so determines a representation $\phi$ of $C^{*}(\Ga,A\otimes B)$.  Then on the range $P$ of $\Phi$, $\pi\circ \Phi^{-1}=\phi$, and it follows that 
$\norm{\Phi^{-1}(z)}\leq\norm{z}$ for all $z\in P$.  So $\Phi$ is isometric, and so extends to an isomorphism from 
$A\otimes C^{*}(\Ga,B)$ onto $C^{*}(\Ga,A\otimes B)$, giving (\ref{eq:equal}).  For the continuity of 
$\bf{F}$, we take $A=C(I)$ where $I$ is some $[a,b]$.  Then using (\ref{eq:equal}), let 
$k\in {\bf F}(IB)=C^{*}(\Ga,IB)\cong C(I,C^{*}(\Ga,B))$.  When $k$ belongs (under the isomorphism $\Phi$) to the dense subalgebra 
$C(I)\otimes_{alg} r_{c}^{*}B$ of $C^{*}(\Ga,IB)$, we use (\ref{eq:phitild}) to show that 
$\hat{k}(t_{0})=\widetilde{ev_{t_{0}}}(k)=k(t_{0})$.   By uniform convergence in $C(I,C^{*}(\Ga,B))$, the same is true for 
$k\in C^{*}(\Ga,IB)$, and ${\bf F}$ is continuous.

(b) For exactness, one modifies the proof by N. C. Phillips of the corresponding result for the group case 
(\cite[Lemma 2.8.2]{PhillipsLNS}).  Let 
\[  0\to J\overset{\chi}{\to}  A \overset{\phi}{\to} B\to 0              \]
be a $\Ga$-equivariant short exact sequence of $\Ga$-algebras.  With $j=\bf{F}(\chi), \psi=\bf{F}(\phi)$, we have to show that 
\[  0\to C^{*}(\Ga,J)\overset{j}{\to}  C^{*}(\Ga,A) \overset{\psi}{\to} C^{*}(\Ga,B)\to 0              \]
is a short exact sequence of $\css$.  So we have to show that (1) $j$ is injective, (2) $\psi\circ j=0$, 
(3) $\ker \psi\subset j(C^{*}(\Ga,J))
$, and (4)
$\psi$ is surjective.

(1) Let $\Phi$ be a representation of $C^{*}(\Ga,J)$.  Let $(L,\pi)$ be a disintegration of $\Phi$ on 
$H=L^{2}(X,\mu;\mfH)$ as earlier.  Regarding the elements of $A$ as multipliers on $J$ in the obvious way, $\pi$ exends to homomorphism $\pi'$ of $A$.  Further, 
\begin{equation}   \label{eq:bai}
\pi'(a)\xi=\lim \pi(ae_{n})\xi 
\end{equation}
where $\{e_{n}\}$ is sequence that is a bounded approximate identity for $J$.  It follows that every $\pi'(a)$ is decomposable, and there is, for each $x\in X$, a representation $\pi'_{x}$ of $A_{x}$ on $H_{x}$ 
such that for each $a\in A$,
\[   \pi'(a)=\int^{\oplus} \pi'_{x}(a_{x})\,d\mu(x).                         \]
Further, for a. e. $x$, $\pi'_{x}$ is non-degenerate.  Next, the restriction of $\pi'$ to $J$ is just $\pi$ 
and so by the uniqueness a. e. of the decomposition of a decomposable operator (\cite[II, 2, 3, Corollary]{D1}) and after
removing a null set from $U$,  we can suppose that $\pi'_{u}$ restricts to $\pi_{u}$ for all $u\in U$.  Then 
$\pi_{u}'(a_{u})=\lim \pi_{u}(a_{u}(e_{n})_{u})$ in the strong operator topology.  We claim that the 
$\{\pi'_{u}\}$ are covariant for the $L_{\ga}$'s.  Indeed, let $a\in A$.  Then with convergence in the strong operator topology,\\ 
$L_{\ga}\pi'(a_{s(\ga)})L_{\ga^{-1}}=\lim L_{\ga}\pi_{s(\ga)}(a_{s(\ga)}(e_{n})_{s(\ga)})L_{\ga^{-1}}
=\lim \pi_{r(\ga)}(\al_{\ga}(a_{s(\ga)}(e_{n})_{s(\ga)}))$\\
$=\lim \pi_{r(\ga)}(\al_{\ga}(a_{s(\ga)}))(\al_{\ga}((e_{n})_{s(\ga)}))
=\pi'_{r(\ga)}(\al_{\ga}(a_{s(\ga)})).$
(Here we use the fact that $\{\al_{\ga}((e_{n})_{s(\ga)})\}$ is a bounded approximate identity for $J_{r(\ga)}$.)  So the pair $(L,\pi')$ is a covariant representation of $(\Ga, A)$ and its integrated form $\Phi'$ is a representation of 
$C^{*}(\Ga,A)$.  Further, since $\chi$ is the identity map, $\Phi(g)=\Phi'(j(g))$ for all 
$g\in r_{c}^{*}J\subset r_{c}^{*}A$.  It follows that 
$\norm{g}\leq \norm{j(g)}$ for all $g\in r_{c}^{*}J$, and by the continuity of $j$, this inequality extends to 
$C^{*}(\Ga,J)$, and $j$ is injective.

(2) $\phi\circ \chi=0$, ${\bf F}(0)=0$ and ${\bf F}$ is a functor.

(3) From (1), $j$ identifies $C^{*}(\Ga,J)$ with a closed ideal of $C^{*}(\Ga,A)$.  Let $g_{0}\in C^{*}(\Ga,A)$ and suppose that $g_{0}\notin C^{*}(\Ga,J)$.  Then there exists a representation $\Phi$ of $C^{*}(\Ga,A)$ such that $\Phi$ annihilates $C^{*}(\Ga,J)$ while 
$\Phi(g_{0})\ne 0$.  Let $(L,\pi)$ be the covariant representation of $(\Ga,A)$ associated with $\Phi$.  
If $h\in C_{c}(X,J^{\sharp})$, $F\in r_{c}^{*}A$, then $hF\in r_{c}^{*}J$, and so by (\ref{eq:mult'}), $\pi(h)=0$.  
So $\pi_{\mid J}=0$ and $\pi$ determines a $C_{0}(X)$-representation $\pi_{1}$ of $A/J=B$.  Also for 
$a\in A$, since $\pi_{1}\circ \phi=\pi$,
\begin{equation}   \label{eq:phif} 
(\pi_{1})_{x}(\phi_{x}(a_{x}))=\pi_{x}(a_{x}),
\end{equation}
It is easy to check that the pair $(L,\pi_{1})$ is covariant for $(\Ga,B)$.  Let $\Phi_{1}$ be the representation of 
$C^{*}(\Ga,B)$ that is the integrated form of $(L,\pi_{1})$.  A simple calculation using 
(\ref{eq:integform}) and (\ref{eq:phif}) gives $\Phi=\Phi_{1}\circ \psi$.  Since $\Phi(g_{0})\ne 0$, we must have 
$\psi(g_{0})\ne 0$.  So $\ker \psi\subset C^{*}(\Ga,J)$.  

(4) Let $F\in r_{c}^{*}B$, $\eps>0$.  Then there exist $F_{i}\in C_{c}(\Ga), b_{i}\in B$ such that 
$\norm{\sum_{i} b_{i}\otimes_{C_{0}(X)} F_{i} - F}<\eps$ in $r_{c}^{*}B$.  By multiplying by a fixed function $g\in C_{c}(\Ga)$ with 
$g=1$ on the support $C$ of $F$, we can suppose that 
there is a fixed compact set $K$, independent of $\eps$, containing the supports of $F$ and the 
$F_{i}$'s.  Then $\norm{\sum_{i} b_{i}\otimes_{C_{0}(X)} F_{i} - F}_{I} \leq \eps\sup_{x}\la^{x}(K\cup K^{-1})$.
Since $\sum_{i} b_{i}\otimes_{C_{0}(X)} F_{i}\in \psi(r_{c}^{*}A)$, it follows that $\psi$ is surjective.
\end{proof}

For later use, in the argument of (a) above, we can take in (\ref{eq:equal}) $A=C_{0}(T)$ to obtain that 
for $k\in C_{0}(T,B)$, the function 
\begin{equation}  \label{eq:c0cont}
\hat{k}\in C_{0}(T,{\bf F}(B)).
\end{equation}


\section{The descent homomorphism}

The theory of $\Ga$-equivariant asymptotic homomorphisms was developed by 
R. Popescu (\cite{Popescu}).  (The case where $\Ga$ is a locally compact group was treated in \cite{GHT}.)
Recall first that in the non-equivariant case, one is given two $\css$ $A, B$.  Let $T=[1,\infty)$.  One defines  
$\mathfrak{A}B=C_{b}(T,B)/C_{0}(T,B)$.  
(The algebras $\mathfrak{A}^{n}B$ ($n\geq 2$) are defined inductively: 
$\mathfrak{A}^{n}B=\mathfrak{A}(\mathfrak{A}^{(n-1)}B)$, but for convenience, we restrict our discussion to the case 
$n=1$.)  An asymptotic morphism is a $*$-homomorphism $\phi$ from $A$ into $\mathfrak{A}B$.   

The theory of asymptotic morphisms in the $C_{0}(X)$-category requires natural and simple modifications 
(\cite{ParkT,Popescu}).  The algebras $A, B$ are, of course, taken to be $C_{0}(X)$-algebras.  However
$C_{b}(T,B)$ is not a $C_{0}(X)$-algebra under the natural homomorphism $\theta:C_{0}(X)\to ZM(C_{b}(T,B))$, where  
$(\theta(f)F)(t)=fF(t)$ ($f\in C_{0}(X), F\in C_{b}(T,B)$).   The reason is that $C_{0}(X)C_{b}(T,B)\ne C_{b}(T,B)$.
Instead, one replaces $C_{b}(T,B)$ by its submodule 
$C_{b}^{X}(T,B)=C_{0}(X)C_{b}(T,B)$ which is a $C_{0}(X)$-algebra.  To ease the notation, we write $C_{b}(T,B)$
instead of $C_{0}(X)C_{b}(T,B)$ when no misunderstanding can arise.  Recall (earlier) that 
$C_{0}(T,B)$ is always a $C_{0}(X)$-algebra with $(C_{0}(T,B))_{x}=C_{0}(T, B_{x})$.  One defines 
$\mathfrak{A}_{X}B$, which, abusing notation slightly, will be abbreviated to $\mfA B$, to be the quotient
$C_{b}(T,B)/C_{0}(T,B)$; then $\mathfrak{A}B$ is a $C_{0}(X)$-algebra.  A 
{\em $C_{0}(X)$-asymptotic morphism} is defined to be a $C_{0}(X)$-morphism $\phi:A\to \mathfrak{A}B$.  

Now suppose that $A, B$ are $\Ga$-algebras.  We would like $C_{b}(T,B), \mathfrak{A}B$ to be $\Ga$-algebras in a natural way 
so that we can define $\Ga$-equivariant asymptotic morphisms from $A$ to 
$\mathfrak{A}B$.  As we will see, there is a technical difficulty in defining the appropriate $\Ga$-actions, and indeed, even in the group case of \cite{GHT}, continuous versions of $C_{b}(T,B), \mathfrak{A}B$ are required.  The $\css$
$(C_{b}(T,B))_{x}$ make sense, of course, since $C_{b}(T,B)$ is now a $C_{0}(X)$-algebra.  The problem is to obtain a natural $\Ga$-action on $C_{b}(T,B)$: how does one define the 
$\al_{\ga}:(C_{b}(T,B))_{s(\ga)}\to (C_{b}(T,B))_{r(\ga)}$?   To deal with this it is natural to try to replace 
$C_{b}(T,B)_{x}$ by $C_{b}(T,B_{x})$ and $C_{b}(T,B)^{\sharp}$ by the bundle 
$C_{b}(T,B^{\sharp})=\sqcup_{x\in X} C_{b}(T,B_{x})$; for, using the $\Ga$-action on $B$, 
$C_{b}(T,B^{\sharp})$ is a $\Ga$-space of $\css$ in the natural way: 
\begin{equation}   \label{eq:raction}
\al_{\ga}(h_{s(\ga)})(t)= \al_{\ga}(h_{s(\ga)}(t)) 
\end{equation}
where, of course, $h_{s(\ga)}\in C_{b}(T,B_{s(\ga)})$.  For each $x$, there is a canonical homomorphism 
$R_{x}:C_{b}(T,B)\to C_{b}(T,B_{x})$ given by: $R_{x}(F)(t)=F(t)_{x}$.  Note that $R_{x}(fF)=f(x)R_{x}F$
($f\in C_{0}(X)$).  Since  $R_{x}(I_{x}C_{b}(T,B))=0$, the map $R_{x}$
descends to a $^{*}$-homomorphism, also denoted $R_{x}$ from $C_{b}(T,B)_{x}$ into $C_{b}(T,B_{x})$. Since 
$R_{x}(C_{0}(T,B))\subset C_{0}(T,B_{x})$, it also induces a $^{*}$-homomorphism, 
$i_{x}:\mathfrak{A}B\to \mathfrak{A}(B_{x})$.  Then $i_{x}$ determines a $^{*}$-homomorphism, also denoted
$i_{x}:(\mathfrak{A}B)_{x}\to \mathfrak{A}(B_{x})$: $i_{x}(\ov{F}_{x})=\ov{R_{x}F}$ where 
$\ov{F}=F+C_{0}(T,B)$ and for $g\in C_{b}(T,B_{x})$, we set $\ov{g}=g+C_{0}(T,B_{x})$.  If we knew that 
$R_{x}, i_{x}$ were onto isomorphisms, then we could identify $C_{b}(T,B)_{x}$ with $C_{b}(T,B_{x})$ and 
$(\mathfrak{A}B)_{x}$ with $\mathfrak{A}(B_{x})$ and be able to define (at least algebraically) actions of $\Ga$ on 
$C_{b}(T,B), \mathfrak{A}B$. 

Unfortunately, the $\cs$ $C_{b}(T,B)$ is too big for this to work (as we will see below).  However, there is a very useful, simple relation between sections of the bundles $C_{b}(T,B)^{\sharp}$ and $C_{b}(T,B^{\sharp})$ which we now describe.  For each 
$F\in C_{b}(T,B)=C_{0}(X,C_{b}(T,B)^{\sharp})$ (Theorem~\ref{th:asharp}), define a section $RF$ of the bundle 
$C_{b}(T,B^{\sharp})$ by setting: $RF(x)=R_{x}F=R_{x}F_{x}$.   Let $S_{0}(X,C_{b}(T,B^{\sharp}))$ be the $\cs$ of sections
of the bundle $C_{b}(T,B^{\sharp})$ that vanish at infinity.  The support $supp\,H$ of a section
$H\in S_{0}(X,C_{b}(T,B^{\sharp}))$ is the closure in $X$ of the set $\{x\in X: H(x)\ne 0\}$.  For 
$F\in C_{b}(T,B)$, the support $supp\,F$ of $F$ is the ($X-$) support of the section $x\to F_{x}$ (not the support of $F$ as a function of $t$).

\begin{proposition}   \label{prop:Rinjec}
The map $R:C_{b}(T,B)\to S_{0}(X,C_{b}(T,B^{\sharp}))$ is a support preserving $^{*}$-isomorphism.
\end{proposition}
\begin{proof}
The only non-trivial part of the proof that $R$ is a $^{*}$-isomorphism is to show $R$ is one-to-one.  Suppose then that
$RF=0$.  Then for fixed $t$, $F(t)_{x}=0$ in $B_{x}$ for all $x$.  So $F(t)=0$ by Theorem~\ref{th:asharp}, and so 
$F=0$.  Now let $F\in C_{b}(T,B)$ be general.  Since $R_{x}$ is norm-decreasing, we obtain that 
$supp\,RF\subset supp\,F$.  Now suppose that $supp\,RF\ne supp\,F$.  Then there exists an open set $V$ in $X$
such that $V\cap supp\,RF=\emp$ and $f\in C_{c}(V)$ such that $fF\ne 0$.  Then $R(fF)=f(RF)=0$ and we contradict the 
$^{*}$-isomorphism property.
\end{proof}

In most of what follows, we will replace $C_{b}(T,B)$ by a smaller $C_{0}(X)$-subalgebra $\mfB$ that contains
$C_{0}(T,B)$.  The constructions above for $C_{b}(T,B)$ go through for $\mfB$.  Let $\mfB_{T}=\sqcup_{x\in X} R_{x}\mfB$, a bundle of $\css$ that is a subbundle of $C_{b}(T,B^{\sharp})$.   As above, we obtain, for each $x$, a $^{*}$-homomorphism $R_{x}:\mfB\to R_{x}\mfB$, which descends to $R_{x}:\mfB_{x}\to R_{x}\mfB$.
Then $\mfB/C_{0}(T,B)$ is a $C_{0}(X)$-algebra, with $(\mfB/C_{0}(T,B))_{x}=\mfB_{x}/C_{0}(T,B_{x})$.
We obtain a canonical homomorphism $i_{x}:(\mfB/C_{0}(T,B))_{x}\to R_{x}(\mfB_{x})/C_{0}(T,B_{x})$.
We write $\mfA_{\mfB} B, \mfA_{\mfB} B_{x}$ in place of $\mfB/C_{0}(T,B), R_{x}(\mfB_{x})/C_{0}(T,B_{x})$.  Note that $R_{x}, i_{x}$ are onto, but unfortunately, are not usually injective.

For example, consider the case where $X=[0,1]$ and $B=C([0,1])$.  Then the function $F$ on $T$, where $F(t)(x)=\sin(tx)$, belongs to $C_{b}(T,B)$ using the mean-value theorem, and $R_{0}(F)(t)=\sin(t0)=0$.  So $R_{0}F=0$.  But $F_{0}\ne 0$.   For otherwise, $F=fF'$ for some $f\in I_{0}, F'\in C_{b}(T,B)$, and so 
$\sup_{t\in T}\left|sin(tx)\right|\to 0$ as $x\to 0$.  
But if $x$ is not zero, then $\sup_{t\in T}\left|sin(tx)\right|=1$.  It is obvious that the image $\ov{F}_{0}$ of 
$F$ in $(\mathfrak{A} B)_{0}$ is non-zero yet $i_{0}(\ov{F}_{0})=0$ in $\mathfrak{A}B_{0}$.  So $i_{0}$ is not 
injective as well.

We now look at the question: when are the $R_{x}$'s $^{*}$-isomorphisms for $\mfB$ as above?  If the latter holds, then 
every $R_{x}$ must be an isometry on $\mfB_{x}$, and it follows that the map 
$x\to \norm{R_{x}F}$ is upper semicontinuous for all $F\in \mfB$.  Here is the converse.  

\begin{proposition}   \label{prop:injec}
Suppose that the map $y\to \norm{R_{y}F}$ is upper semicontinuous for all $F\in \mfB$.  Then 
for every $y\in X$, $R_{y}$ is a $^{*}$-isomorphism on $\mfB_{y}$.
\end{proposition}
\begin{proof}
Let $F\in \mfB$, $y\in X$ and 
$\eps>0$.  Since $\norm{R_{x}F}\leq \norm{F_{x}}\to 0$ as $x\to \infty$, there exists a compact subset $C$ of $X$ 
such that $\norm{R_{x}F}<\eps$ for all $x\in X\setminus C$. Suppose that $R_{y}F=0$.
By upper semicontinuity, there exists an open neighborhood 
$U$ of $y$ such that $\norm{R_{x}F}<\eps$ for all $x\in U$.  We can suppose that $U\subset C$.  Let 
$h\in C_{c}(X)$ be such that $0\leq h\leq 1$, $h(y)=0$ and $h(x)=1$ for all $x\in C\setminus U$.  Also 
$\norm{R_{x}(F - hF)}=0$ if $x\in C\setminus U$, $\leq\norm{R_{x}F}<\eps$ if $x\in U$, $< \eps$ if 
$x\in X\setminus C$.  So $\norm{R_{x}(F - hF)}<\eps$ for all $x\in X$. 
Note next that $\norm{R_{x}(F - hF)}=(1-h(x))\norm{R_{x}F}$.  So for each $t$, 
$\norm{F(t)_{x} - h(x)F(t)_{x}}<\eps$ for all $x$, and so by Theorem~\ref{th:asharp}, $\norm{F(t)-hF(t)}\leq\eps$.  
So $\norm{F - hF}\leq \eps$.  Since $h\in I_{y}$, $\norm{F_{y}}\leq \eps$, and since $\eps$ was arbitrary, 
$F_{y}=0$.  So 
$R_{y}:\mfB_{y}\to C_{b}(T,B_{y})$ is injective and so isometric.  
\end{proof}

We note that under the condition of Proposition~\ref{prop:injec}, every $i_{y}$ is also isometric.
Now let $B$ be a $\Ga$-algebra.   We have the canonical action of $\Ga$ on the bundle $\mfB_{T}$: 
$\al_{\ga}(F_{s(\ga)})(t)=\al_{\ga}(F(t)_{s(\ga)})$, $F\in \mfB$ provided that $\mfB$ is {\em $\Ga$-invariant}, i.e.  
$\al_{\ga}(R_{s(\ga)}\mfB)=R_{r(\ga)}\mfB$ for all $\ga$. In that situation, we look for an action $\bt$ of $\Ga$ on $\mfB$ -  not necessarily continuous - related to $\al$ on 
$\mfB_{T}$; precisely, we want $\bt$ to {\em cover} $\al$ in the sense that for each $\ga$, the following diagram commutes:
\[
\begin{CD}
\mfB_{s(\ga)}                   @>\bt_{\ga}>>            \mfB_{r(\ga}  \\
@VR_{s(\ga)}VV                                            @VVR_{r(\ga)}V\\
R_{s(\ga)}\mfB_{s(\ga)}         @>>\al_{\ga}>          R_{r(\ga)}\mfB_{r(\ga)}   
\end{CD}
\]
i.e., 
\begin{equation}  \label{eq:covering}
R_{r(\ga)}\bt_{\ga}=\al_{\ga}R_{s(\ga)}   
\end{equation}  
on $\mfB_{s(\ga)}$.  If this condition is satisfied and if $\bt$ restricted to $C_{0}(T,B)$ is the canonical 
$\Ga$-algebra action on $C_{0}(T,B)$ ((\ref{eq:gactionc0})), then we say that $\mfB$ has a {\em covering $\Ga$-action}.  In the case of a continuous covering $\Ga$-action, we can extend the continuity condition for a continuous functor ${\bf F}$ from finite intervals $I=[a,b]$ to the infinite interval $T$.  For such an $I$, the canonical $\Ga$-action on $IB$ 
(also as in (\ref{eq:gactionc0})) will be denoted by $\ga\to \al_{\ga}'$.  It is determined by: 
$R_{r(\ga)}\al_{\ga}'=\al_{\ga}R_{s(\ga)}$, where $\al_{\ga}'(h_{s(\ga)})(t)=\al_{\ga}(h(t)_{s(\ga)})$ (using $R_{x}$ to identify $(IB)_{x}$ with $IB_{x}$).

\begin{proposition}   \label{prop:Tcont}
Let $\mfB\subset C_{b}(T,B)$ have a covering $\Ga$-action and $\bf{F}$ be a continuous functor as in the third section. Then for all $k\in \bf{F}(\mfB)$, the function $\hat{k}\in C_{b}(T,\bf{F}(B))$, and if $\mfB=C_{0}(T,B)$, then 
$\hat{k}\in C_{0}(T,\bf{F}(B))$.
\end{proposition}
\begin{proof}
Let $I$ be a closed bounded interval inside $T$ and $\rho:\mfB\to IB$ be the restriction map.  We show first that 
$\rho$ is a $\Ga$-homomorphism.  Note first that $\rho$ is a $C_{0}(X)$-homomorphism: for since 
$fF(t)=f\cdot F(t)$ ($f\in C_{0}(X), F\in \mfB$) for all $t$, we get $\rho(fF)=f\rho(F)$.  Also, for $t\in I$, 
$R_{x}\rho_{x}F_{x}(t)=F(t)_{x}=R_{x}F_{x}(t)$.  Using the above and (\ref{eq:covering}), for $t\in I$,
\[ R_{r(\ga)}\rho_{r(\ga)}\bt_{\ga}F_{s(\ga)}(t)=R_{r(\ga)}\bt_{\ga}F_{s(\ga)}(t)=\al_{\ga}[R_{s(\ga)}F_{s(\ga)}(t)]\]
\[  =\al_{\ga}[R_{s(\ga)}\rho_{s(\ga)}F_{s(\ga)}](t) = R_{r(\ga)}\al_{\ga}'\rho_{s(\ga)}F_{s(\ga)}(t).  \]
Since $R_{r(\ga)}$ is a bijection, we obtain $\rho_{r(\ga)}\bt_{\ga}= \al_{\ga}'\rho_{s(\ga)}$, so that $\rho$ is a $\Ga$-homomorphism.  With 
$I$ a single point $t_{0}$, we get that $ev_{t_{0}}$ is also a $\Ga$-homomorphism.  Now let 
$I=[a,b]\subset T$ and $\rho$ be as above.  For $t\in I$, $ev_{t}=ev_{t}\circ\rho$, and so for 
$k\in {\bf F}(\mfB)$, $\hat{k}(t)={\bf F}(ev_{t})(k)={\bf F}(ev_{t})({\bf F}(\rho)(k))$ which is continuous in $t$ by the continuity of $\bf{F}$.  So $\hat{k}\in C_{b}(T,{\bf F}(B))$.  The last part is just (\ref{eq:c0cont}).
\end{proof}

One natural way in which a covering $\Ga$-action on some $\mfB$ can arise is when, for each $x$, $R_{x}$ is an isomorphism on 
$\mfB_{x}$ and $\mfB_{T}$ is $\Ga$-invariant, in which case we can identify the bundles $\mfB^{\sharp}$ and 
$\mfB_{T}$  and obtain a $\Ga$-action on $\mfB$: effectively $\bt=\al$ in this case.  An example of this is the situation of 
Proposition~\ref{prop:unifcont} below.  However, there are many cases where the $R_{x}$'s are not isomorphisms but we can still find a covering action on some reasonable 
$\mathfrak{B}$.  For instance, in the example above with $X=[0,1]$, where we take $\Ga=X$ (a groupoid of units), we can take $\mfB=C[0,1]$ and $\bt_{x}(F_{x})=F_{x}$!  Of course, this example in trivial, but as we will see in 
Theorem~\ref{th:thfinal}, such a $\mfB$ always exists under very general circumstances.  

Suppose now that $\Ga$ has local $r-G$-sets, and that $\mfB$ has a covering $\Ga$-action.  By 
Proposition~\ref{prop:gacs}, $\mfB^{cont}$ is a $\Ga$-$\cs$.   By definition, the $\bt$ action restricts to give the canonical 
$\Ga$-action on $C_{0}(T,B)$ which is continuous.  So $C_{0}(T,B)\subset \mfB^{cont}$.  So we can define  
$\mfA_{\mfB^{cont}} B=\mfB^{cont}/C_{0}(T,B)$.  If $A$ is also a $\Ga$-$\cs$, then an {\em equivariant asymptotic morphism} from $A$ to $B$ (relative to $\mfB^{cont}$) is just a $\Ga$-homomorphism 
$\phi:A\to \mfA_{\mfB^{cont}} B$.  In the locally compact group case, one takes $\mfB=C_{b}(T,B)$ and the asymptotic algebra $\mfA_{\mfB^{cont}} B=C_{b}(T,B)^{cont}/C_{0}(T,B)$ is the same as the $\mfA B$ of 
\cite[p.7]{GHT}.  In that case, there is a descent {\em functor} for $\Ga$-$\css$ using as morphisms 
homotopy classes of $\Ga$-homomorphisms into the asymptotic algebra (\cite[Theorem 4.12]{GHT}).  However, since, for completely general 
$\Ga$,  we do not have available a canonical $\mfB^{cont}$, it does not make sense to talk of ``the 
functor $B\to \mfA_{\mfB^{cont}} B$''.   (However, for a very wide class of groupoids $\Ga$, we do obtain a canonical asymptotic algebra and a functor from Theorem~\ref{th:thfinal} below - it would be interesting to know if the complete theory of the descent functor can be developed for this class of groupoids as in \cite{GHT}.)  Instead at present, we avoid a functorial description of the descent functor and give a direct, weaker version of the descent homomorphism which is adequate for a number of applications.

\begin{theorem} \label{th:descentcovering}
Suppose that $\Ga$ has local $r-G$-sets, and that $\mfB$ is a $C_{0}(X)$ subalgebra of $C_{b}(T,B)$ containing 
$C_{0}(T,B)$ and with a covering $\Ga$-action.  Let $A$ be 
a $\Ga$-$\cs$ and $\phi:A\to \mfA_{\mfB^{cont}} B$ be a $\Ga$-homomorphism.  Then there exists a canonical descent homomorphism (dependent on $\mfB$) given by: 
$i\circ {\bf F}(\phi):C^{*}(\Ga,A)\to \mfA C^{*}(\Ga,B)$.
\end{theorem}
\begin{proof}
The proof is effectively the same as for the group case (\cite[Theorem 4.12]{GHT}).  Let ${\bf F}$ be as in
Theorem~\ref{th:exact}.  From the exactness of ${\bf F}$, we get the short exact sequence:
\[ 0\to {\bf F}(C_{0}(T,B))\to {\bf F}(\mfB^{cont})\to {\bf F}(\mfA_{\mfB^{cont}} B)\to 0.  \]
We also have the short exact sequence for the ordinary $\cs$ ${\bf F}(B)$:
\[ 0\to C_{0}(T,{\bf F}(B))\to C_{b}(T,{\bf F}(B))\to \mathfrak{A} {\bf F}(B)\to 0.         \]
Continuity (Proposition~\ref{prop:Tcont}) gives $^{*}$-homomorphisms $i_{0}:{\bf F}(C_{0}(T,B))\to C_{0}(T,{\bf F}(B))$, 
$i_{b}:{\bf F}(\mfB^{cont})\to C_{b}(T,{\bf F}(B))$ with $i_{0}$ the restriction of $i_{b}$ to 
${\bf F}(C_{0}(T,B))$, and these induce a $^{*}$-homomorphism 
$i:{\bf F}(\mfA_{\mfB^{cont}} B)\to \mathfrak{A}{\bf F}(B)$.  Next we have a $^{*}$-homomorphism 
${\bf F}(\phi): {\bf F}(A)\to {\bf F}(\mfA_{\mfB^{cont}} B)$.  So
$i\circ {\bf F}(\phi):{\bf F}(A)\to \mathfrak{A}{\bf F}(B)$ is a $^{*}$-homomorphism.
\end{proof}  

Before discussing our main theorem giving a canonical $\mfB$ we look at a situation in which there is a very simple 
$\mfB$ available.  For motivation, Proposition~\ref{prop:injec} suggests that we should look for a $\mfB$ with elements 
$F$ for which the map $x\to \norm{R_{x}F}$ is upper semicontinuous.  I have been unable to find a canonical such $\mfB$ in general.  However, in cases that arise in practice - in particular, when $\Ga$ is  discrete, or when $\Ga$ is a locally compact group (the case of \cite{GHT}) or when $B=\mathcal{C}(E)=C_{0}(\R)\hat{\otimes}C_{0}(E,Cliff(E))$ 
($E$ a $\Ga$-vector bundle) - there is a natural such $\mfB$ available.  (The last of these cases is needed for the groupoid version of the infinite dimensional Bott periodicity theorem of Higson, Kasparov and Trout (\cite{HKT}).)  Intuitively, we wish to exclude from $\mfB$ functions such as $\sin(xt)$ by requiring {\em uniformly continuity} in the 
$X$-direction.  This requires a strong condition on $B$, but the groupoid $\Ga$ has to satisfy only the weak local 
$r-G$-set condition of Proposition~\ref{prop:gacs}.

We know that $B$ is isomorphic to the $\cs$ $C_{0}(X,B^{\sharp})$ of the $C^{*}$-bundle 
$(B^{\sharp},p)$.  We assume now that this $C^{*}$-bundle $B^{\sharp}$ is a locally trivial $C^{*}$-bundle with ($\cs$) fiber $C$.  Precisely, a chart $(U,\eta)$ is given by an open subset $U$ of $X$ together with a fiber preserving homeomorphism $\eta$ from $p^{-1}(U)$ onto $U\x C$ with each $\eta_{x}$ a $^{*}$-isomorphism from $B_{x}$ onto $C$ 
($x\in U$).  Local triviality means that the chart sets $U$ cover $X$.  (In particular, no structure group condition is imposed.)  For such a chart $(U,\eta)$ and $F\in C_{b}(T,B)$, $R_{x}F\in C_{b}(T,B_{x})$ and so for $x\in U$,
$\eta_{x}\circ R_{x}F\in C_{b}(T,C)$.  
We say that $F$ is {\em uniformly continuous} (for $X$) if the map $\eta R(F): x\to \eta_{x}\circ R_{x}F$ is continuous from 
$U$ to $C_{b}(T,C)$ for every chart $(U,\eta)$.   
It is easily checked that the set $\mfB=UC_{b}(T,B)$ of uniformly continuous functions is a $C^{*}$-algebra, and a 
$C_{0}(X)$-subalgebra of $C_{b}(T,B)$.   Note that the equality $R_{x}\mfB=C_{b}(T,B_{x})$ below shows that 
this $\mfB$ is ``big''.

\begin{proposition}    \label{prop:unifcont}
Let $\mfB=UC_{b}(T,B)$.  Then $C_{0}(T,B)\subset \mfB$, and for $F\in \mfB$, the map $x\to \norm{R_{x}F}$ is continuous.
The maps $R_{x}, i_{x}$ are $^{*}$-isomorphisms.   Further $R_{x}(\mfB)=C_{b}(T,B_{x})$, and $\mfB$ is $\Ga$-invariant (so that trivially $\mfB$ has a covering $\Ga$-action).  
\end{proposition}
\begin{proof} Let 
$k\in C_{0}(T,B)$.  To show that $k\in UC_{b}(T,B)$, we can suppose that $k=h\otimes b$ where $h\in C_{c}(T)$ and 
$b\in C_{c}(X)B$.  Using a partition of unity, we can suppose that 
there is a chart $(U,\eta)$ and a compact subset $K$ of $U$ such that $k_{x}=0$ for all $x\in X\setminus K$.  Then 
$\eta R(k)\in C_{c}(T\x U,C)$, and $k\in \mfB$.  So $C_{0}(T,B)\subset \mfB$.  Now let $F\in \mfB$.  Then on $U$,
the map $x\to \norm{\eta_{x}\circ R_{x}F}=\norm{R_{x}F}$ is continuous.  By Proposition~\ref{prop:injec}, each 
$R_{x}, i_{x}$ is an isomorphism.  By definition, $R_{x}(\mfB)\subset C_{b}(T,B_{x})$.   To show equality, we just have to show that if $H\in C_{c}(U,C_{b}(T,C))$ then there exists $F\in C_{b}(T,B)$ such that $\eta R(F)=H$ (so that 
$F\in \mfB$).    For then we can take any $\phi\in C_{b}(T,B_{x})$, take $g\in C_{c}(U)$ with $g(x)=1$ and 
$H(y)=g(y)\eta_{x}\circ \phi$ to get $R_{x}F=\phi$.   

To show that such an $F$ exists, fix $t$.  The map $y\to H_{y}(t)\in C$ is continuous on $U$, and so 
$y\to \eta_{y}^{-1}H_{y}(t)$ is a continuous section of $B^{\sharp}$ supported on $U$.  By Theorem~\ref{th:asharp}, 
there exists $b_{t}\in B$ such that $(b_{t})_{y}=\eta_{y}^{-1}H_{y}(t)$.  Define $F(t)=b_{t}$.  Then $F$ is bounded since 
$\norm{(b_{t})_{y}}\leq \norm{H_{y}(t)}\leq \norm{H}_{\infty}$.  Last $F$ is continuous.  Indeed, 
given $\eps>0$, there exists a $\de>0$ such that for all $y\in U$, $\norm{H_{y}(t)-H_{y}(s)}<\eps$ whenever
$\left| t-s \right|<\de$.  Let $\left| t-s \right|<\de$.  Then for all $y$, $\norm{(b_{t})_{y} - (b_{s})_{y}}<\eps$. 
Now take the supremum over $y$ to get $F$ continuous at $t$.
\end{proof}

It follows from Theorem~\ref{th:descentcovering} that if $\Ga$ has local $r-G$-sets, then 
$UC_{b}(T,B)$ determines a descent homomorphism: this can be regarded as the natural descent homomorphism for 
such a special $B$.   We now show that under a very mild condition on $\Ga$ and with $B$ completely general, 
there {\em always} exists a canonical $\mfB$ with a covering $\Ga$-action giving a 
descent homomorphism.  The algebra $\mfB$ is functorial.  The condition that we need on $\Ga$ is that it have {\em local
$G$-sets}, i.e. local transversals for $r$ and $s$ simultaneously.  The algebra $\mfB$ consists (roughly) of those functions with a transversally continuous action, and going to $\mfB^{cont}$ then gives continuity of the action in every direction.  We now make these ideas precise.

A subset $W$ of $\Ga$ is called a {\em G-set} (cf. \cite[p.10]{rg})
if the restrictions $r_{W}, s_{W}$ of $r, s$ to $W$ are homeomorphisms onto open subsets of $X$.  The family of 
$G$-sets in $\Ga$ is denoted by $\mathcal{G}$.

\begin{proposition}   \label{prop:G-sets}
$\mathcal{G}$ is an inverse semigroup under pointwise product and inversion.
\end{proposition}
\begin{proof}
The discrete case is given in \cite[Proposition 2.2.3]{Paterson}.  For the topological conditions, we need to show that 
for $W, W_{1}, W_{2}\in \mathcal{G}$, the bijections $r_{W_{1}W_{2}}, s_{W_{1}W_{2}}, r_{W^{-1}},\\ 
s_{W^{-1}}$ are homeomorphism onto open subsets of $X$.  We will prove this for $r_{W_{1}W_{2}}$ leaving the rest to the reader.  First, 
$r(W_{1}W_{2})=r_{W_{1}}s_{W_{1}}^{-1}(r(W_{2})\cap s(W_{1}))$ which is open in $X$.  Next let $x_{n}\to x$ in 
$r(W_{1}W_{2})$.  We write uniquely $x_{n}=r(w_{1}^{n}w_{2}^{n}), x=r(w_{1}w_{2})$ where $w_{1}^{n}, w_{1}\in W_{1}, 
w_{2}^{n}, w_{2}\in W_{2}$.  Then $(r_{W_{1}})^{-1}(x_{n})=w_{1}^{n}\to (r_{W_{1}})^{-1}(x)\\
=w_{1}$, and similarly, 
$w_{2}^{n}\to w_{2}$.   So $w_{1}^{n}w_{2}^{n}\to w_{1}w_{2}$ and $r_{W_{1}W_{2}}^{-1}$ is continuous.
\end{proof}

We note that if $W\in \mathcal{G}$, the map $r_{W}s_{W}^{-1}$ is a homeomorphism from $s(W)$ onto $r(W)$, and 
there is a nice formula for $r_{W}s_{W}^{-1}$: 
\[   r_{W}s_{W}^{-1}(x)=WxW^{-1}.                      \]
We will say that $\Ga$ has {\em local $G$-sets} if $\cup\mathcal{G}=\Ga$, i.e. every element of $\Ga$ belongs to a $G$-set.  This property is satisfied by most groupoids that arise in practice (e.g. r-discrete groupoids, transformation group groupoids, many (all?) Lie groupoids).  For motivation for the following, suppose that $W\in \mathcal{G}$.  Suppose that $\mfB\subset C_{b}(T,B)$ has a covering $\Ga$-action $\bt$ that makes it into a $\Ga$-algebra.  For $F\in \mfB$, 
 the map $r(\ga)\to \bt_{\ga}(F_{s(\ga)})$ is continuous and so will come from an element $F'$ of 
$C_{b}(T,B)$ (at least after multiplying by a function in $C_{c}(s(W))$).   Now from Proposition~\ref{prop:Rinjec}, we can recover 
$F'$ from $RF'$ and at the $R$-level, we do have the action $\al_{\ga}$.  The idea then is to consider functions $F$ for which there is an $F'$ that goes down under $R$ to the function $r(\ga)\to \al_{\ga}(R_{s(\ga)}F)$ and define 
$\bt_{\ga}(F_{s(\ga)})=F'_{r(\ga)}$.  We also insist that this definition is independent of the choice of $W$.
This does not ensure a continuous action but only one continuous along the
$G$-sets.  However, the algebra of such functions is the largest subalgebra $C_{b}^{\Ga}(T,B)$ of $C_{b}(T,B)$ admitting a reasonable covering action that {\em is} continuous along the $G$-sets.  We now develop the theory of $C_{b}^{\Ga}(T,B)$.

We identify $C_{b}(T,B)$ with 
$C_{0}(X,(C_{b}(T,B))^{\sharp})$.  Abbreviate $C_{b}(T,B)$ to $C_{b}$, and let $F\in C_{b}$.  
We say that $F\in C_{c}^{G}=C_{c}^{G}(\Ga,B)$ if:
\be
\item $F\in C_{c}(X)C_{b}$;
\item for all $W\in \mathcal{G}$ with $supp\,F\subset s(W)$, there exists $F^{W}\in C_{b}$ such that
\[     RF^{W}=R_{W}F              \]
where $R_{W}F(x)=\al_{\ga}R_{s(\ga)}F$ if $x=r(\ga)$ for some $\ga\in W$, and is $0$ otherwise;
\item (uniqueness) if $W_{1}, W_{2}\in \mathcal{G}$, $\ga_{0}\in W_{1}\cap W_{2}$ and 
$supp\,F\subset s(W_{1})\cap s(W_{2})$, then 
\[           (F^{W_{1}})_{r(\ga_{0})}=(F^{W_{2}})_{r(\ga_{0})}.   \]
\ee
We say that $F\in C_{b}^{\Ga}=C_{b}^{\Ga}(T,B)$ if $fF\in C_{c}^{G}$ for all $f\in C_{c}(X)$.  

\begin{theorem}   \label{th:thfinal}
$C_{b}^{\Ga}(T,B)$ is a $C_{0}(X)$-subalgebra of $C_{b}(T,B)$ with a covering $\Ga$-action given by:
\begin{equation}  \label{eq:btga0}
\bt_{\ga_{0}}F_{s(\ga_{0})}=(fF)^{W}_{r(\ga_{0})}    
\end{equation}
where $W\in \mathcal{G}, \ga_{0}\in W$ and $f\in C_{c}(s(W))$ is such that $f(s(\ga_{0}))=1$.  
\end{theorem}
\begin{proof}
We prove the theorem in five stages.

(a) {\em Let $F\in C_{c}^{G}$ and $W\in \mathcal{G}$ with $supp\,F\subset s(W)$.  Then $F^{W}$ is unique, and 
$supp\,F^{W}=W (supp\,F) W^{-1}$.  Also $\norm{F^{W}}=\norm{F}$.}\\
The uniqueness of $F^{W}$ follows from Proposition~\ref{prop:Rinjec}.  The same proposition gives that 
$supp\,F=supp\,RF, supp\,R_{W}F=supp\,F^{W}$.  Next, since $\al_{\ga}$ is an isometry, for $\ga\in W$,
$R_{s(\ga)}F\ne 0$ if and only if $R_{W}F(r(\ga))\ne 0$, so that $R_{x}F\ne 0$ if and only if 
$R_{W}F(WxW^{-1})\ne 0$ for $x\in s(W)$.  Closing up gives $supp\,R_{W}F=W (supp\,RF) W^{-1}$.  For the last part, use Proposition~\ref{prop:Rinjec} and the fact that the $\al_{\ga}$'s are $^{*}$-isomorphisms. 

(b) {\em Let $F\in C_{c}^{G}$ and $W\in \mathcal{G}$ with $supp\,F\subset s(W)$.  Then $F^{W}\in C_{c}^{G}$, and for 
$V\in \mathcal{G}$ with $supp\,F^{W}\subset s(V)$, $(F^{W})^{V}=F^{VW}$.  Further, $F^{*}\in C_{c}^{G}$ and $(F^{*})^{W}=(F^{W})^{*}$.}\\
By (a), $supp\,F^{W}=W(supp\,F)W^{-1}\subset V^{-1}V$, and so $supp\,F\subset s(VW)$.  Now
$R_{V}F^{W}(x)=\al_{\ga_{1}}R_{s(\ga_{1})}F^{W}$ if $x=r(\ga_{1})$ for some $\ga_{1}\in V$, and is $0$ otherwise.  Next
$\al_{\ga_{1}}R_{s(\ga_{1})}F^{W}=\al_{\ga_{1}}(\al_{\ga_{2}}R_{s(\ga_{2})}F)$ if $s(\ga_{1})=r(\ga_{2})$ for some (unique)
$\ga_{2}\in W$ and is $0$ otherwise.  Since 
$\al_{\ga_{1}}(\al_{\ga_{2}}R_{s(\ga_{2})}F)=\al_{\ga_{1}\ga_{2}}R_{s(\ga_{1}\ga_{2})}F$ and 
$r(\ga_{1})=r(\ga_{1}\ga_{2})$, it follows that $R_{V}F^{W}=RF^{VW}$.  So $(F^{W})^{V}=F^{VW}$.  The uniqueness condition (3) of the definition of $C_{c}^{G}$ with respect to $F^{W}$ follows from the corresponding property for $F$.  The last part of the proof of (b) is easy.

(c) {\em Let $F\in C_{b}^{\Ga}, W\in \mathcal{G}$, $f\in C_{c}(s(W))$ and $supp\,F\subset s(W)$ be compact.  Then 
$(fF)^{W}=f^{W}F^{W}$ and belongs to $C_{b}^{\Ga}$, and
$R_{W}(fF)=R(f^{W}F^{W})$ where $f^{W}\in C_{c}(r(W))$ is given by: $f^{W}(y)=f(W^{-1}yW)$.}\\
Note that by definition, $fF\in C_{b}^{G}$.  Next, if $\ga\in W$ and $y=r(\ga)$, then 
$R_{W}(fF)(y)=\al_{\ga}(R_{s(\ga)}(fF))=f(s(\ga))R_{W}(F)(y)=f^{W}(r(\ga))R_{W}F(r(\ga))
=R_{y}(f^{W}F^{W})$.  Both $R_{W}(fF), R(f^{W}F^{W})$ vanish at $y$ if $y\notin r(W)$, and so $R_{W}(fF)=R(f^{W}F^{W})$.
By (a), $(fF)^{W}=f^{W}F^{W}$.  Next
we show that $(fF)^{W}\in C_{b}^{\Ga}$.  For let $g\in C_{c}(X)$ and 
$h\in C_{c}(r(W))$ be such that $h=1$ on $supp\,F^{W}=W (supp\,F) W^{-1}$, a compact subset of $r(W)$ (using (a)).
Then $g(fF)^{W}=(gh)(fF)^{W}=((gh)^{W^{-1}})^{W}(fF)^{W}=((gh)^{W^{-1}}fF)^{W}\in C_{c}^{G}$ by (b).  So
$(fF)^{W}\in C_{b}^{\Ga}$. 

(d) {\em $C_{b}^{\Ga}$ is a $C_{0}(X)$-subalgebra of $C_{b}$ that contains $C_{0}(T,B)$.}\\
Trivially, $0\in C_{b}^{\Ga}$.  Let $F_{1}, F_{2}\in C_{b}^{\Ga}, f\in C_{c}(X), W\in \mathcal{G}$ with
$C=supp\,(f(F_{1}+F_{2}))\subset s(W)$.  Choose $f'\in C_{c}(s(W))$ so that $f'=1$ on $C$.  Then 
$f(F_{1}+F_{2})=f'fF_{1}+f'fF_{2}$ with $supp\,f'fF_{1}\cup supp\,f'fF_{2}\subset s(W)$.  Since 
$F_{1}, F_{2}\in C_{b}^{\Ga}$, we get $f'fF_{1}, f'fF_{2}\in C_{c}^{G}$.  Then 
$R_{W}(fF_{1}+fF_{2})=R_{W}(f'fF_{1}) + R_{W}(f'fF_{2})=R((f'fF_{1})^{W}+(f'fF_{2})^{W})$ and we obtain
$(fF_{1}+fF_{2})^{W}=(f'fF_{1})^{W}+(f'fF_{2})^{W}$.  For uniqueness, let $W'\in \mathcal{G}$, $C\subset s(W')$
and $\ga_{0}\in W\cap W'$.  Then we can choose $f'\in C_{c}(s(W)\cap s(W'))$ and obtain uniqueness using 
$((f'fF_{1})^{W})_{r(\ga_{0})}=((f'fF_{1})^{W'})_{r(\ga_{0})},
((f'fF_{2})^{W})_{r(\ga_{0})}=((f'fF_{2})^{W'})_{r(\ga_{0})}$.  So $F_{1}+F_{2}\in C_{b}^{\Ga}$.  Next, 
$F_{1}F_{2}\in C_{b}^{\Ga}$.  The proof is similar to that for the sum above.  Let $W\in \mathcal{G}$ with 
$C=supp\,fF_{1}F_{2}\subset s(W)$.  Choose $f'\in C_{c}(s(W))$ with $f'=1$ on $C$.  Then 
$fF_{1}F_{2}=(f'fF_{1})(f'F_{2})$ and we can use $R_{W}(fF_{1}F_{2})=R_{W}(f'fF_{1})R_{W}(f'F_{2})$ to get
$(fF_{1}F_{2})^{W}=(f'fF_{1})^{W}(f'F_{2})^{W}$.  The remaining details are left to the reader, as is also the proof that 
$F_{1}^{*}\in C_{b}^{\Ga}$ (use (b)).  So $C_{b}^{\Ga}$ is a $^{*}$-subalgebra of $C_{b}$. 

Next we show that $C_{b}^{\Ga}$ is complete.  Let $F_{n}\to F$ in $C_{b}$ with every $F_{n}\in C_{b}^{\Ga}$.  We show that 
$F\in C_{b}^{\Ga}$.  Let $f\in C_{c}(X)$ and $W\in \mathcal{G}$ be such that $supp\,fF\subset s(W)$.  We can suppose that 
$f\in C_{c}(s(W))$.  Let $g\in C_{c}(s(W))$ be $1$ on the support of $f$.  Using (a) and Proposition~\ref{prop:Rinjec}, 
$\norm{gF_{n}-gF_{m}}=\norm{R(gF_{n}-gF_{m})}=\norm{(gF_{n})^{W} - (gF_{m})^{W}}\to 0$ as $n, m\to\infty$.  So there exists $F''\in C_{b}$ such that $(gF_{n})^{W}\to F''$ in norm.  By the continuity of $R_{W}, R$, we get for $\ga\in W$,
$R_{W}(fF)(r(\ga))\\
=\lim R_{W}(f(gF_{n}))(r(\ga))=\lim R_{r(\ga)}(f^{W}(gF_{n})^{W})=R_{r(\ga)}(f^{W}F'')$.  We take
$(fF)^{W}\\
=f^{W}F''$.  Uniqueness is proved using uniqueness for the $F_{n}$'s and a simple limit argument.  To prove that 
$C_{0}(T,B)\subset C_{b}^{\Ga}$, let $F\in C_{0}(T,B)$, $f\in C_{c}(X)$, $W\in \mathcal{G}$, $supp(fF)\subset s(W)$.
Then we take $(fF)^{W}(r(\ga))=\al_{\ga}(fF)_{s(\ga)}$ ($\ga\in W$), where $\al$ is the canonical $\Ga$-algebra action.
Uniqueness is obvious.  So $fF\in C_{c}^{G}$ and $F\in C_{b}^{\Ga}$.

Last we have to show that $C_{b}^{\Ga}$ is a $C_{0}(X)$-algebra.  Let $F\in C_{b}^{\Ga}$ and $h\in C_{0}(X)$.  Trivially, if $f\in C_{c}(X)$ then $f(hF)=(fh)F\in C_{c}^{G}$.  So $hF\in C_{b}^{\Ga}$.  To prove that the action of $C_{0}(X)$ on 
$C_{b}^{\Ga}$ is non-degenerate use the fact that $C_{b}=C_{0}(X)C_{b}$.

(e) {\em The $\bt_{\ga}$'s give a covering $\Ga$-action on $C_{b}^{\Ga}$.}\\
We first show that, for given $F\in C_{b}^{\Ga}$, the right hand side of (\ref{eq:btga0}) is well-defined.  
Let $f,W$ be as in (\ref{eq:btga0}) and let $W'\in \mathcal{G}$, $\ga_{0}\in W'$ and $f'\in C_{c}(s(W'))$ be such that 
$f'(s(\ga_{0}))=1$.  Suppose first that $W=W'$.  Let $g\in C_{c}(s(W))$ be such that $g=1$ on 
$supp\,f\,\cup\, supp\,f'$.  Then $fF=fgF$, $f'F=f'gF$, and using (c), 
$(fF)^{W}_{r(\ga_{0})}=f^{W}(r(\ga_{0}))(gF)^{W}_{r(\ga_{0})}=(gF)^{W}_{r(\ga_{0})}=(f'F)^{W}_{r(\ga_{0})}$.  For the case 
$W\ne W'$, find $h\in C_{c}(s(W)\cap s(W'))$ such that $h(s(\ga_{0}))=1$.  Then using uniqueness, 
$(fF)^{W}_{r(\ga_{0})}=(hF)^{W}_{r(\ga_{0})}=(hF)^{W'}_{r(\ga_{0})}=(f'F)^{W}_{r(\ga_{0})}$.  Next we have to show that
the right hand side of (\ref{eq:btga0}) depends only on the coset $F_{s(\ga_{0})}$.  So let $F'=F+gF_{1}$ where
$g\in I_{s(\ga_{0})}$ and $F_{1}\in C_{b}^{\Ga}$.  Then 
$(fF')^{W}_{r(\ga_{0})}=(fF)^{W}_{r(\ga_{0})}+g^{W}(r(\ga_{0})(fF_{1})_{r(\ga_{0})}=(fF)^{W}_{r(\ga_{0})}$.  We now show that $\ga\to \bt_{\ga}$ defines an algebraic action on $C_{b}^{\Ga}$.  It is simple, using (b) and the proof of (d) to show that each $\bt_{\ga}$ is a $^{*}$-homomorphism.  

Next we show that 
$\bt_{\ga_{0}}\bt_{\ga_{1}}=\bt_{\ga_{0}\ga_{1}}$ whenever $s(\ga_{0})=r(\ga_{1})$.  Let $F\in C_{b}^{\Ga}$ and 
let $W\in \mathcal{G}$ contain $\ga_{1}$ and $f\in C_{c}(s(W))$ be such that $f(s(\ga_{1}))=1$.  Then 
$\bt_{\ga_{1}}F_{s(\ga_{1})}=(fF)^{W}_{r(\ga_{1})}$.  By (c), $(fF)^{W}\in C_{b}^{\Ga}$ so that 
$(fF)^{W}_{r(\ga_{1})}\in (C_{b}^{\Ga})_{r(\ga_{1})}$.  Now let $V\in \mathcal{G}$ be such that $\ga_{0}\in V$.  Since 
$s(\ga_{0})=r(\ga_{1})$ belongs to $s(V)\cap r(W)$, we can find
$g\in C_{c}(s(V)\cap r(W))$ such that $g(s(\ga_{0}))=1$.  Then 
$\bt_{\ga_{0}}(\bt_{\ga_{1}}F_{s(\ga_{1})})=(g(fF)^{W})^{V}_{r(\ga_{0})}
=(((g^{W^{-1}}fF)^{W})^{V})_{r(\ga_{0})}=((g^{W^{-1}}fF)^{VW})_{r(\ga_{0}\ga_{1})}$ using (b).  Noting that 
$(g^{W^{-1}}f)(s(\ga_{0}\ga_{1}))=g(Ws(\ga_{0}\ga_{1})W^{-1})f(s(\ga_{1}))=g(r(\ga_{1}))f(s(\ga_{1}))=1$
and $\ga_{0}\ga_{1}\in VW$, we see that $((g^{W^{-1}}fF)^{VW})_{r(\ga_{0}\ga_{1})}$ is just 
$\bt_{\ga_{0}\ga_{1}}F_{s(\ga_{1})}$ and  $\bt_{\ga_{0}}\bt_{\ga_{1}}=\bt_{\ga_{0}\ga_{1}}$.   If $x\in X$, then trivially (using $W=X$), $\bt_{x}:(C_{b}^{\Ga})_{x}\to (C_{b}^{\Ga})_{x}$ is the identity, and it follows that 
$\ga\to \bt_{\ga}$ is an algebraic $\Ga$-action on $C_{b}^{\Ga}$.  To show that the $\Ga$-action is covering, in an obvious notation, $R_{r(\ga)}(\bt_{\ga}F_{s(\ga)})=R_{r(\ga)}(fF)^{W}=R_{W}(fF)(r(\ga))
=\al_{\ga}(R_{s(\ga)}(fF))=\al_{\ga}R_{s(\ga)}(fF)_{s(\ga)}=\al_{\ga}R_{s(\ga)}F_{s(\ga)}$. It is left to the reader to 
check that $\bt_{\ga}$ restricted to $C_{0}(T,B)$ gives the canonical $\Ga$-action on $C_{0}(T,B)$.
\end{proof}

Since $\Ga$ has local $G$-sets, it satisfies the condition of Proposition~\ref{prop:gacs}, and so 
$C_{b}^{\Ga}(T,B)^{cont}$ is a $\Ga$-algebra.  This can be regarded as the canonical $\Ga$-algebra $\mfB$ for groupoids with local $G$-sets.  Indeed, the map $B\to C_{b}^{\Ga}(T,B)^{cont}$ is functorial in the category of $\Ga$-algebras with $\Ga$-homomorphisms as morphisms.  To see this, let $B_{1}, B_{2}$ be $\Ga$-algebras and $\phi:B_{1}\to B_{2}$ be a $\Ga$-homomorphism.  Let $\mfB_{1}= C_{b}^{\Ga}(T,B_{1}), \mfB_{2}=C_{b}^{\Ga}(T,B_{2})$.  Define 
$\tilde{\phi}:C_{b}(T,B_{1})\to C_{b}(T,B_{2})$ by: $\tilde{\phi}(F)(t)=\phi(F(t))$.  Then $\psi=\tilde{\phi}$ is a 
$C_{0}(X)$-homomorphism.  One readily checks that $R_{x}(\psi(F))=\phi_{x}\circ R_{x}F$.  Let $F\in \mfB_{1}$.
Let $W\in \mathcal{G}$, $f\in C_{c}(X)$ and $supp(fF)\subset s(W)$.  Then for $\ga\in W$,
$\al_{\ga}R_{s(\ga)}(f\psi(F))=\al_{\ga}R_{s(\ga)}(\psi(fF))=\al_{\ga}(\phi_{s(\ga)}\circ R_{s(\ga)}(fF))
=\phi_{r(\ga)}\circ (\al_{\ga}R_{s(\ga)}(fF))=\phi_{r(\ga)}\circ R_{r(\ga)}((fF)^{W})=R_{r(\ga)}(\psi((fF)^{W}))$.
So $(f\psi(F))^{W}=\psi((fF)^{W})$ exists, and uniqueness is easily checked from that for $F$.  So 
$f\psi(F)\in C_{c}^{G}(T,B_{2})$ and 
$\psi(F)\in C_{b}^{\Ga}(T,B_{2})$.  Next, if $\ga_{0}\in s(W)$ and $f\in C_{c}(s(W))$ with $f(s(\ga_{0}))=1$, we get
$\bt_{\ga_{0}}(\psi_{s(\ga_{0})}F_{s(\ga_{0})})=\bt_{\ga_{0}}((\psi(F))_{s(\ga_{0})})
=(f\psi(F))^{W}_{r(\ga_{0})}=(\psi((fF)^{W}))_{r(\ga_{0})}=\psi_{r(\ga_{0})}(fF)^{W}_{r(\ga_{0})}=
\psi_{r(\ga_{0})}\bt_{\ga_{0}}F_{s(\ga_{0})}$, so that $\psi$ is equivariant from $\mfB_{1}$ to $\mfB_{2}$.
Last, using equivariance, $\psi:\mfB_{1}^{cont}\to \mfB_{2}^{cont}$ and the functorial property follows.

\end{document}